\documentclass[11pt]{article}

\usepackage{amsfonts,amsmath,amsthm,amscd,amssymb,latexsym,amsbsy}

    \usepackage[T2A]{fontenc}
    \usepackage[cp1251]{inputenc}
   \usepackage[russian,english]{babel}

 \usepackage{srcltx}
 \makeatletter
  \renewcommand{\@oddhead}{\hfil\thepage}
  \renewcommand{\@oddfoot}{}
\makeatother

\begin{document}


\setcounter{page}{395}

\renewcommand{\C}{\mathbb{C}} 
\newcommand{\N}{\mathbb{N}}
\newcommand{\R}{\mathbb{R}}
\newcommand{\Z}{\mathbb{Z}}

\newcommand{\sign}{\mathop{\rm sign}\nolimits}
\renewcommand{\Re}{\mathop{\rm Re}\nolimits}
\renewcommand{\Im}{\mathop{\rm Im}\nolimits}


\theoremstyle{plain}
\newtheorem{z_theorem}{Theorem}[section]
\newtheorem{z_lemma}{Lemma}[section]
\newtheorem{z_proposition}{Proposition}[section]
\newtheorem{z_corollary}{Corollary}[section]
\newtheorem{z_definition}{Definition}[section]
\newtheorem{oldtheorem}{Theorem}
\renewcommand{\theoldtheorem}{\Alph{oldtheorem}}

\theoremstyle{definition}
\newtheorem{z_remark}{Remark}[section]

\newtheorem{z_example}{Example}[section]
\newcommand{\keywords}{\textbf{Key words and phrases. }\medskip}
\newcommand{\subjclass}{\textbf{2000 MSC. }\medskip}
\renewcommand{\abstract}{\textbf{Abstract. }\medskip}
\renewcommand{\refname}{\textbf{References}}
\numberwithin{equation}{section}

\title{On Exponential Type Entire Functions without Zeros in the Open Lower Half-plane}


\author{Viktor P. Zastavnyi}

\date{({\it Presented by L. A. Pastur})\vspace*{.5cm}\\ {\footnotesize Version published in Ukrainian Mathematical Bulletin, vol. 3(2006), \No~3, 395--422}}


 \maketitle
\noindent\begin{abstract}
  We obtain sufficient conditions for an exponential type entire
  function not to have zeros in the open lower half-plane. An exact
  inequality containing the real and imaginary parts of such
  functions and their derivatives restricted to the real axis is
  deduced.
\end{abstract}

\vspace*{.5cm}\noindent\subjclass{30C15, 30D15, 42A82, 60E10.}

 \vspace*{.5cm}\noindent\keywords{Entire function, Hermite--Biller theorem, positive definite function, Fourier transform.}


\thispagestyle{empty}

\section{Introduction. Formulation of the Main Results}

In this paper, we study entire functions that do not have zeros in
the open lower half-plane, $\Im z<0$. The results applied to
algebraic polynomials give the known Hermite--Biller theorem. An
extension of this theorem to arbitrary entire functions was carried
out in works of M.~G.~Krein, B.~Ya.~Levin, and N.~N.~Meiman, for
more details see, for example,~\cite{Lev1,Lev2,Ostr, Cheb_Me}.

An entire function $f$ is called an exponential type entire function
(in~\cite{Lev1}, such functions are called finite exponent
functions), if there exist numbers $A>0,\ B>0$ such that $|f(z)|\le
A e^{B|z|} $ holds for all $z\in\C$. The exact lower bound of such
numbers $B$ is denoted by $\sigma (f)\ge 0$ and is called a type of
the function $f$. Denote by $E_{\sigma},\ \sigma \ge 0,$ the class
of entire functions $f$ of exponential type $\sigma (f)\le \sigma$.

The purpose of this work is to prove Theorems~\ref{z:t1},
\ref{z:t2}, \ref{z:t3}, and~\ref{z:t4}.

\pagebreak
\begin{z_theorem}\label{z:t1}
  Let the following conditions hold: $\rm 1)$
  $\omega(z)\!=\!P(z)+iQ(z)$, where $P,Q$ are real\footnote{We call
    a function real if it takes real values on the real axis.}
  functions of the class $E_{\sigma}$, $\sigma > 0$,
  and $\omega(x)=o(x)$, $x\to\pm\infty$ on the real axis; $\rm 2)$
  the following inequality holds for some $\tau\in\R$:
  \begin{equation}\label{6}
    E(x):= P(x)\cos(\sigma x+\tau)+Q(x)\sin(\sigma  x+\tau)\ge
    0,\qquad x\in\R.
  \end{equation}
  Let $d(x):=P(x)Q'(x)-P'(x)Q(x)$. Then the following holds.
  \begin{itemize}
  \item[{\bf 1)}] For all $x\in\R$,
    \begin{multline}\label{7}
      4\sigma d(x)\ge
      \bigl\{(\sigma P(x)+Q'(x))\sin(\sigma x+\tau)\\+(P'(x)-\sigma Q(x))\cos(\sigma
      x+\tau)\bigr\}^2.
    \end{multline}

  \item[{\bf 2)}] The following conditions are equivalent:
    \begin{itemize}  \item[{\rm i)}]
      inequality~\eqref{7} becomes equality;

    \item[{\rm ii)}] for some $c\ge 0$, $\beta\in\R$, we have the
      identity $E(x)\equiv c\sin^2(\sigma x+\tau+\beta)$;

    \item[{\rm iii)}] for some $\beta\in\R$ and all $k\in\Z$, we
      have $E\bigl(\frac{k\pi-\beta-\tau}{\sigma}\bigr)=0$;

    \item[{\rm iv)}] for some $c\ge 0$ and $\beta,\gamma\in\R$, we
      have
      \begin{equation}\label{8}
        \begin{aligned}
          P(x)&\equiv c\sin\beta \sin(\sigma x+\tau+\beta)+\gamma\sin(\sigma  x+\tau),\\
          Q(x)&\equiv c\cos\beta \sin(\sigma x+\tau+\beta)-\gamma\cos(\sigma
          x+\tau).
        \end{aligned}
      \end{equation}
      In this case, $d(x)\equiv \gamma^2\sigma$.
    \end{itemize}

  \item[{\bf 3)}] If for any $\alpha\in\R$ and $c\ge 0$,
    $E(x)\not\equiv c\sin^2(\sigma x+\tau+\alpha)$, then
    inequality~\eqref{7} becomes equality for some $x=x_0\in\R\iff
    E(x_0)=0$.
  \end{itemize}
\end{z_theorem}

For an entire function $\omega(z)=P(z)+iQ(z)$, where $P(z)$ and
$Q(z)$ are real entire functions, set
$\overline{\omega}(z):=P(z)-iQ(z)$. This is an entire function
obtained from $\omega(z)$ by replacing all the coefficients in its
expansion in the powers $\{z^{n-1}\}_{n\in\N}$ with their complex
conjugates. It is clear that
$\overline{\omega}(z)=\overline{\omega(\overline{z})}$.

\begin{z_definition}\label{z:d1}
  An entire function $\omega(z)$ is called a class $HB$ function if
  it does not have zeros in the closed lower half-plane $\Im z\le
  0$, and $\bigl|\frac{\omega(z)}{{\overline\omega}(z)}\bigr|<1$ if
  $\Im z>0$.
\end{z_definition}
\pagebreak
\begin{z_definition}\label{z:d2}
  An entire function $\omega(z)$ is called a class $\overline{HB}$
  function if it does not have zeros in the open lower half-plane
  $\Im z< 0$, and
  $\bigl|\frac{\omega(z)}{{\overline\omega}(z)}\bigr|\le 1$ if $\Im
  z>0$.
\end{z_definition}

Equality in Definition~\ref{z:d2} can hold only if $\omega(z)$ is a
real function up to a constant factor. Such functions are called
trivial class $\overline{HB}$ functions. It is clear that {\it
  $\omega\in HB\iff$ the function $\omega\in \overline{HB}$ does not
  have real zeros and is not trivial.  }

\begin{z_theorem}\label{z:t2}
  Let conditions of Theorem~{\rm \ref{z:t1}} be satisfied and assume
  that the function $\omega$ is not real up to a constant factor.
  Then we have the following.
  \begin{itemize}
  \item[{\bf 1)}] $\omega\in\overline{HB}$. If the function $\omega$
    has real zeros, then they are simple.

  \item[{\bf 2)}] $d(x_0)= 0$ for some $x_0\in{\mathbb{R}} \iff
    \omega(x_0)=0 $. If a number $x_0\in\R$ is a zero of the
    function $\omega$, then the number $x_0$ is a zero of the
    function $d$ of multiplicity $2$.

  \item[{\bf 3)}] For all $n\in{\mathbb{N}}$, $\omega^{(n)}\in HB$.
  \end{itemize}
\end{z_theorem}

In Corollary~\ref{sl2}, we give examples where conditions of
Theorems~\ref{z:t1} and \ref{z:t2} hold. Note that
Theorems~\ref{z:t1} and~\ref{z:t2} cease to hold if the condition
$\omega(x)=o(x)$, $x\to\pm\infty$, is replaced with the condition
$\omega(x)=O(x)$, $x\to\pm\infty$, see Remark~\ref{re2}.

For a given function $\mu$ that has bounded variation on the segment
$[0,\sigma]$, $\sigma>0$, we will consider the following entire
functions:
\begin{equation}\label{1}
  \begin{gathered}
    F(z):=\int\limits_{0}^{\sigma} e^{izt} \,d\mu(t),\qquad
    G(z):=\int\limits_{0}^{\sigma} \cos{zt} \,d\mu(t),\\
    H(z):=\int\limits_{0}^{\sigma} \sin{zt} \,d\mu(t),
  \end{gathered}
\end{equation}
\begin{equation}\label{2}
\begin{aligned}
\Delta(z)&:=G(z)H'(z)-G'(z)H(z),\\
 h_{\alpha}(z)&:=G(z)\cos{\alpha} - H(z)\sin{\alpha} ,
 \end{aligned}
\end{equation}
\begin{equation}\label{3}
C(z):=-\int\limits_{0}^{\sigma} \cos{zt} \,d\mu(\sigma-t),\qquad
S(z):=-\int\limits_{0}^{\sigma} \sin{zt} \,d\mu(\sigma-t).
\end{equation}
It is clear that the function $\mu$ can always be regarded as left
continuous in every point of the interval $(0,\sigma)$. If
$F(z)\not\equiv ce^{i\alpha z}$, then the function $F$ has
infinitely many zeros, see for example~\cite{Leont}. A function of
type~(\ref{1}) is used for solving many problems in analysis, for
example, in spectral theory, in the theory of
differential-difference equations, in the theory of positive
definite functions, in the Fourier analysis, etc., see for
example~\cite{Leont, Luk, Sed2, Tih, Tr}. The distribution of
zeros of such functions was studied in the works of
Hardy~\cite{Hardy}, Polya~\cite{Pol}, Titchmarsh~\cite{Tit},
Cartwright~\cite{Cart1,
  Cart2}, Sedlitskii~\cite{Sed1, Sed3}, and others. In~\cite{Pol},
the case where the function $\mu$ is absolutely continuous and
$d\mu(t)=g(t)dt$ was studied, where the function $g$ is integrable,
positive, and does not decrease on the interval $(0,\sigma)$. It was
proved there that, in this case, all zeros of the function $F$ lie
in the closed upper half-plane $\Im z\ge 0$, and if $g$ is not
piecewise constant with uniformly distributed nodes, then the
function $F$ does not have real zeros. These results of Polya were
extended and made more precise in the works of the
author~\cite{Zast2004_MZ, Zast2004_MFAT}.

\begin{z_theorem}\label{z:t3}
  Let $\mu$ be a real function with bounded variation on the segment
  $[0,\sigma]$, $C(x)\ge 0$ for all $x\in\R$, and $F(z)\not\equiv
  0$. Then the function $F\in\overline{HB}$ is not trivial. All real
  zeros of the function $F$, if there are any, are simple.
\end{z_theorem}

\begin{z_theorem}\label{z:t4}
  Let $\mu$ be a real function with bounded variation on the segment
  $[0,\sigma]$, $S(x)\ge 0$ for all $x>0$, and $F(z)\not\equiv 0$.
  Then all real zeros distinct from $0$ of the function $F$, if
  there are any, are simple, and if the number $x=0$ is a zero of
  the function $F$, then its multiplicity does not exceed $2$.
  Moreover, in such a case, $\mu(\sigma-0)\ge\mu(0) $ and the
  following holds.
  \begin{itemize}
  \item[{\bf 1)}] If $F(0)\in
    (-\infty,0]\cup[\mu(\sigma-0)-\mu(0),+\infty)$, then
    $F\in\overline{HB}$, and the function is not trivial.

  \item[{\bf 2)}] If $F(0)\in ( 0, \mu(\sigma-0)-\mu(0) )$, then the
    function $F$ has exactly one zero in the lower half-plane
    ${\mathop{\rm Im}}\, z<0$, and it is pure imaginary.
  \end{itemize}
\end{z_theorem}

In Section~5, we give examples where conditions of
Theorems~\ref{z:t3} and~\ref{z:t4} hold. Polya's case is contained
in statement~1 of Theorem~\ref{z:t4} (see Example~5.1 in Section~5).
Note that Theorem~\ref{z:t4} realizes two cases (if the quantity
$\mu(\sigma)$ is changed, then the values of $S(x)$ and
$\mu(\sigma-0)$ do not change, and the value
$F(0)=\mu(\sigma)-\mu(0)$ can be made arbitrary). Also note that
$S(x)\ge 0$ for all $x>0\iff$ $\mu(t)-\mu(0)\equiv f(\sigma-t)$ for
$0\le t<\sigma$, where $f$ is an even function that is positive
definite and continuous on ${\mathbb{R}}$, and equal to zero for
$|t|\ge \sigma$ (Lemma~\ref{le2}). A relation between class
$\overline{HB}$ functions of type~\eqref{1} and positive definite
functions is contained in Proposition~\ref{pr7} (see Section~5).

\section{Auxiliary Propositions}

\subsection{Functions~(\ref{1}), (\ref{2}), and~(\ref{3}).}

We have the following:
\begin{equation}\label{a1}
\begin{aligned}
  F(z)&\equiv G(z)+iH(z),\\
  F(z)e^{-i\sigma z}&\equiv -\int\limits_{0}^{\sigma}
  e^{-izt}\,d\mu(\sigma-t)\equiv C(z)-iS(z),
 \end{aligned}
\end{equation}
\begin{equation}\label{a2}
\begin{aligned}
G(z)&\equiv C(z)\cos{\sigma z}+S(z)\sin{\sigma z},\\
 H(z)&\equiv C(z)\sin{\sigma z}-S(z)\cos{\sigma z},
 \end{aligned}
\end{equation}
\begin{equation}\label{a3}
  G(z)\cos(\sigma z+\tau)+H(z)\sin(\sigma z+\tau)\equiv C(z)\cos{\tau}-S(z)\sin{\tau},
\end{equation}
\begin{equation}\label{a4}
\begin{aligned}
  h_{\alpha}(z)&\equiv C(z)\cos(\sigma z+\alpha)+S(z)\sin(\sigma
  z+\alpha),\\
  h_{\alpha}(z)h'_{\beta}(z)-h'_{\alpha}(z)h_{\beta}(z)&\equiv \Delta(z)\sin(\alpha-\beta),
  \end{aligned}
\end{equation}
\begin{equation}\label{a5}
\Delta(x)\equiv C'(x)S(x)-C(x)S'(x)+\sigma
\left(C^2(x)+S^2(x)\right).
\end{equation}
These identities can be directly obtained from~(\ref{1}), (\ref{2}),
and~(\ref{3}).
\begin{z_lemma}\label{le1}
  {\bf 1)} {\rm i)} $G(z)\equiv 0\iff F(z)\equiv 0$.  {\rm ii)}
  $H(z)\equiv 0\iff F(z)\equiv c$.  {\rm iii)} $C(z)\equiv 0\iff
  F(z)\equiv 0$.  {\rm iv)} $S(z)\equiv 0\iff F(z)\equiv ce^{i\sigma
    z}$.  {\rm v)} $h_{\alpha}(z)\equiv 0$ for some $\alpha\iff
  G(z)\cos\alpha\equiv H(z)\sin\alpha\equiv 0$.

  {\bf 2)} The function $F$ is real up to a constant factor $\iff
  F(z)\equiv c$.

  {\bf 3)} $x^n(C(x)\cos\tau-S(x)\sin\tau)\equiv c\sin^2(\sigma
  x+\tau+\alpha)$ for some $\tau,\alpha,c\in\C$, $n\in\Z\iff
  C(x)\cos\tau - S(x)\sin\tau\equiv{0}\iff C(x)\cos\tau\equiv
  S(x)\sin\tau\equiv{0}$.

  {\bf 4)} If $F(z)\equiv ce^{(i\alpha+\beta)z}$ for some
  $\alpha,\beta\in\R$, $c\ne 0$, then $\beta=0$ and
  $\alpha\in[0,\sigma]$.

  {\bf 5)} If $F(z)\equiv ce^{i\alpha z}$ for some $c\ne 0$,
  $\alpha\in[0,\sigma]$, and for all $x>0$ the inequality
  $C(x)\cos\tau-S(x)\sin\tau\ge 0$ holds, then $\alpha=\sigma$.

  {\bf 6)} If $F(z)\not\equiv 0$ and $C(x)\cos\tau-S(x)\sin\tau\ge
  0$ for $x>0$, then the function $F$ is not real up to a constant
  factor, and $h_{\alpha}(z)\not\equiv 0$.

  {\bf 7)} If for some $a,b,c,d,e\in\R$ and any $t>e$, the
  inequality $f(t):=a\sin^2(t+b)+c\cos t+d\sin t\ge 0$ holds, then
  $c=d=0$ and $a\ge 0$.

  {\bf 8)} If $\mu$ is real and $F(0)\ne 0$, then $F(z)\not\equiv c
  R(z)(z+i\xi)$, where $\xi\in\R$ and $R(z)$ is a real entire
  function.
\end{z_lemma}

\begin{proof}
  Let us prove statement~1).

  i) If $G(z)\equiv 0$, then
  $\int_{0}^{\sigma} t^{2p}d\mu(t)=0$ for all
  $p\in\Z_+:=\N\cup\{0\}$. By Muntz's theorem, the system of
  powers $\{t^{2p}\}_{p\in\Z_+}$ is dense in $C[0,1]$. Hence,
  $\int_{0}^{\sigma} f(t)d\mu(t)=0 $ for any function $f\in C[0,1]$.
  Thus, $F(z)\equiv 0$. Conversely, it follows from the identity
  $F(z)=G(z)+iH(z)\equiv 0$, since $G$ is even and $H$ is odd, that
  $G(z)\equiv H(z)\equiv 0$.

  \sloppy
  ii) If $H(z)\equiv 0$, then $\int_{0}^{\sigma} t^{2p+1}d\mu(t)=0$,
  $p\in\Z_+$. Then $\int_{0}^{\sigma} f(t)\,t\,d\mu(t)=0 $ for
  any function $f\in C[0,1]$. Hence, $F'(z)\equiv 0$, and so
  $F(z)\equiv c$. Conversely, if $F(z)\equiv c$, then $2iH(z)\equiv
  F(z)-F(-z)\equiv 0$.

  \fussy
  Statements~iii) and~iv) follow from i), ii),
  and~(\ref{a1}). Statement v) follows at once from~(\ref{2}) if we
  recall that $G$ is even and $H$ is odd functions.

  \sloppy
  Let us prove statement 2). Let a function $F$ be real up to a
  constant factor. Without loss of generality, we can assume that
  the function $F$ is real. Let $\mu(t)=\mu_1(t)+i\mu_2(t)$, where
  $\mu_1(t) ,\ \mu_2(t)$ are real functions with bounded variation
  on the segment $[0,\sigma]$. Then $\Im(F(x))\equiv
  \int_{0}^{\sigma} \sin xt\,d\mu_1(t)+\int_{0}^{\sigma} \cos
  xt\,d\mu_2(t)\equiv 0$, $x\in\R$, and, hence, $\int_{0}^{\sigma}
  \sin xt\,d\mu_1(t)\equiv\int_{0}^{\sigma} \cos xt\,d\mu_2(t)\equiv
  0$. It follows from statement 1) that $\int_{0}^{\sigma}
  e^{izt}\,d\mu_1(t)\equiv c$ and $ \int_{0}^{\sigma}
  e^{izt}\,d\mu_2(t)\equiv 0$. Hence, $F(z)\equiv c$. The converse
  is clear.

  \fussy
  Now we prove statement 3). If the indicated identity holds, then
  $c=0$; otherwise the left-hand side contains an entire function of
  exponential type $\le \sigma$, and the type of the function in the
  right-hand side is precisely $2\sigma$. So, $C(x)\cos\tau -
  S(x)\sin\tau\equiv 0$, which is equivalent to two identities,
  $C(x)\cos\tau\equiv S(x)\sin\tau\equiv 0$.

  Let us prove 4). If $\beta\ne 0$, then the right-hand side of the
  identity is unbounded on $\R$, whereas the left-hand side is
  bounded.  Hence, $\beta= 0$ and $F(z)\equiv ce^{i\alpha z}$, $c\ne
  0$. If $\alpha>\sigma$, then the function in the right-hand side
  of the identity has type greater than that of the function in the
  left-hand side. If $\alpha<0$, then $F(iy)=\int_{0}^{\sigma}
  e^{-yt}d\mu(t)\equiv ce^{-\alpha y} $. The left-hand side, as
  $y\to +\infty$, is bounded and the right-hand side is not. Hence,
  $\alpha\in[0,\sigma]$.

  Now, let us prove~5). In this case,
  $C(x)\cos\tau-S(x)\sin\tau\equiv c\cos((\sigma-\alpha)x+\tau)\ge
  0$ for $x>0$, where $c\ne 0$, $\alpha\in[0,\sigma]$. If
  $\alpha<\sigma$, then $\exists\,x_0>0 :
  c\cos((\sigma-\alpha)x_0+\tau)\ne 0$. Then, for $x=x_0$ and
  $x=x_0+\frac{\pi}{\sigma-\alpha}$, the left-hand side of the
  inequality takes values of different sign. Thus, $\alpha=\sigma$.

  Let us prove~6). If a function $F$ is a constant, up to a constant
  factor, then $F(z)\equiv c$, and $c\ne 0$, which can not happen
  (see statement~5 with $\alpha=0$). If $h_{\alpha}(z)\equiv 0$, then
  $G(z)\cos\alpha\equiv H(z)\sin\alpha\equiv 0$ and, hence,
  $F(z)\equiv c$, $c\ne 0$, which is a contradiction.

  Now we prove~7). Let $t_k:=-b+k\pi$, $k\in\Z$. Then, for all
  $k>k_0$, we have $f(t_k)=(-1)^k(c\cos b-d\sin b)\ge 0$ and, hence,
  $c\cos b-d\sin b=0$. Thus, $f(t_k)=0$ for all $k>k_0$, and so,
  $f'(t_k)=(-1)^k(c\sin b+d\cos b)=0$. Consequently, $c=d=0$, which
  means that $a\ge 0$.

  Let us prove~8). Without loss of generality, we can assume that
  the function $\mu$ is left continuous in every point of the
  interval $(0,\sigma)$. Assume that $F(z)\equiv c R(z)(z+i\xi)$,
  where $\xi\in\R$ and $R(z)$ is a real entire function. Since
  $F(0)=ic\xi R(0)\in\R\setminus\{0\}$, we can assume that $c=-i$
  and $\xi\ne 0$. Then $\xi R(z)=\int_{0}^{\sigma} \cos zt\,d\mu(t)$
  and $z R(z)=-\int_{0}^{\sigma} \sin zt\,d\mu(t)=z\int_{0}^{\sigma}
  \cos zt(\mu(t)-\mu(\sigma))\,dt=z\int_{0}^{\sigma} \cos
  zt\,d\mu_1(t)$, where
  $\mu_1(t)=\int_{\sigma}^{t}(\mu(u)-\mu(\sigma))\,du$. Hence, if
  $t\in [0,\sigma]$, we have $\mu(t)-\mu(\sigma)\equiv
  \xi\int_{\sigma}^{t}(\mu(u)-\mu(\sigma))\,du$ which implies that
  $\mu\in C^1[0,\sigma]$ and $\mu(t)-\mu(\sigma)\equiv c_1e^{\xi
    t}$. Consequently, $\mu(t)-\mu(\sigma)\equiv 0$ and $F(z)\equiv
  0$, which contradicts the condition.
\end{proof}

A function $f:\R \to \C$ is called positive definite on $\R$ if, for
any $ n\in \N$, $ \{x_k\}_{k=1}^n \subset \R$, and $ \{c_k\}_{k=1}^n
\subset \C$, the inequality $ \sum_{k,j=1}^n c_k \bar{c}_j
f(x_k-x_j) \ge 0 $ holds. For such functions, $|f(x)|\le f(0)$,
$x\in\R$, and continuity at zero is equivalent to continuity on
$\R$.  By Bokhner--Hinchin theorem, a function $f$ is positive
definite and continuous on $\R$ if and only if
$f(x)=\int_{-\infty}^{+\infty} e^{-iux} \,d\nu (u),$ where $\nu$ is
a nonnegative, finite Borel measure on $\R$. If $f\in C(\R) \cap
L(\R),$ then positive definiteness of a function $f$ is equivalent
to its Fourier transform being nonnegative, that is,
$\widehat{f}(x):=\int_{-\infty}^{+\infty} f(u)e^{-iux}\,du \ge 0$,
$x\in \R $, and in this case, $\widehat{f}\in L(\R)$, see~\cite[Ch.
I, \S 1, Corollary~1.26]{StW}).

\begin{z_lemma}\label{le2}
  Let $\mu$ be a real function of bounded variation on the segment
  $[0,\sigma]$, left continuous in every point of the interval
  $(0,\sigma)$. Then the following holds.
  \begin{itemize}
  \item[{\bf 1)}] $S(x)\ge 0$ for all $x>0\iff$ $\mu(t)-\mu(0)\equiv
    f(\sigma-t)$ for $0\le t<\sigma$, where $f$ is an even function
    that is positive definite and continuous on $\R$, $f(t)=0$,
    $|t|\ge \sigma$. In this case, $f(0)=\mu(\sigma-0)-\mu(0)\ge 0$
    and $\mu(\sigma-0)-\mu(0)=0\iff S(x)\equiv 0\iff F(z)\equiv
    ce^{i\sigma z}, c\in\R$.

  \item[{\bf 2)}] If $S(x)\ge 0$ for all $x>0$ and $F(0)\le 0$, then
    $H'(0)\le 0$. In this case, $F'(0)=0\iff H'(0)=0\iff F(0)=0$ and
    $\int_{0}^{\sigma}f(t)dt= 0$, where $f$ is the corresponding
    function in statement~1).

  \item[{\bf 3)}] If $S(x)\ge 0$ for all $x>0$ and $F(0)\ge
    \mu(\sigma-0)-\mu(0)$, then $H'(0)\ge 0$. In this case,
    $F'(0)=0\iff H'(0)=0\iff F(z)\equiv 0$.
   \end{itemize}
\end{z_lemma}

\begin{proof}
  \sloppy Let us prove statement~1). It follows from~(\ref{3}) and
  the integration by parts formula that $S(x)=xK(x)$, where
  $K(x):=\int_{0}^{\sigma}\cos tx\,(\mu(\sigma-t)-\mu(0))\,dt$. Let
  us first assume that $S(x)\ge 0$ for all $x>0$. The function
  $2K(x)$ is a Fourier transform of a finite function, integrable on
  $\R$, $\mu((\sigma-|t|)_+)-\mu(0)$, which is bounded in a
  neighborhood of zero. Since $K(x)\ge 0$ for all $x\in\R$, we have
  $K\in L(\R)$, see for example~\cite{Tr,Ah}, and hence, for almost
  all $t\in\R$, one can apply the inverse transform formula, see for
  example~\cite{Tr,StW,Ah},
  $$
  \mu((\sigma-|t|)_+)-\mu(0)=\frac{1}{\pi}
  \int\limits_{-\infty}^{+\infty} e^{itx} K(x)\,dx=:f(t)\;.
  $$
  The function $f$ in the right-hand side of the latter identity
  is even, continuous, and positive definite on $\R$. The left-hand
  side equals $0$ for $|t|\ge \sigma$. Hence, continuity of the
  function $f$ implies that $f(t)=0$ for all $|t|\ge \sigma$. Since
  the function $\mu$ is left continuous in every point of the
  interval $(0,\sigma)$, the above identity implies that the
  function $\mu$ is continuous in every point of this interval, and
  $\mu(0+0)-\mu(0)=f(\sigma)=0$. Hence, $\mu\in C[0,\sigma)$ and
  $\mu(t)-\mu(0)\equiv f(\sigma-t)$ for $0\le t<\sigma$. The
  converse is obvious. The first part of statement~1) is proved. The
  second part follows from the inequality $|f(x)|\le f(0)$,
  $x\in\R$, and Lemma~\ref{le1} (statement~1).

  \fussy
  Statement 2) follows at once from the identities $F'(0)=iH'(0)$,
  $F(0)=\mu(\sigma)-\mu(0)$,
  $H'(0)=\int_{0}^{\sigma}t\,d\mu(t)=\sigma
  F(0)-\int_{0}^{\sigma}f(\sigma-t)\,dt$, and the inequality
  $\int_{0}^{\sigma}f(t)\,dt\ge 0$.

  Statement 3) immediately follows from the inequalities $|f(t)|\le
  f(0)$ and $H'(0)=\int_{0}^{\sigma}(F(0)-f(t))\,dt\ge
  \int_{0}^{\sigma}(f(0)-f(t))\,dt\ge 0$. If $H'(0)=0$, then
  $F(0)=f(0)\equiv f(t)$ for $t\in[0,\sigma]$ and, hence, $f(0)=0$,
  $S(x)\equiv 0$, $F(z)\equiv ce^{i\sigma z}$, where $c=F(0)=0$.
\end{proof}

\begin{z_lemma}\label{le3}
  Let $f,g$ be functions positive definite on $\R$ from the class $
  C(\R)\cap L(\R)$, and $f(x)\not\equiv 0$, $g(x)\not\equiv 0$ on
  $\R$. If the function $f$ is finite, then {\bf 1)}
  $\int_{-\infty}^{+\infty}g(x)f(x)\,dx>0$; {\bf 2)} for all
  $\alpha>0$ and $\beta\in\R$,  $|\beta|\le\alpha$, the following
  inequality holds:
  $\int_{-\infty}^{+\infty}e^{-\alpha|x|}(1-\beta|x|)f(x)\,dx>0$.
\end{z_lemma}

\begin{proof}
  Since $\widehat{f}(t)\ge 0$ and $\widehat{g}(t)\ge 0$ for all
  $t\in\R$, it follows that $\widehat{f},\widehat{g}\in L(\R)\cap
  C(\R)$. The multiplication formula gives that
  $$
  \int\limits_{-\infty}^{+\infty}g(x)f(x)\,dx
  =\frac{1}{2\pi}\int\limits_{-\infty}^{+\infty}\widehat{g}(-t)\widehat{f}(t)\,dt\ge
  0\;.
  $$
  If the integral equals $0$, then
  $\widehat{g}(-t)\widehat{f}(t)\equiv 0$ on $\R$. Since
  $g(x)\not\equiv 0$, we have $\widehat{g}(-t)\ne 0$ on some
  interval $(a,b)$, $a<b$. Hence, $\widehat{f}(t)=0$ on $(a,b)$,
  and, since $\widehat{f}$ is entire, $\widehat{f}(t)\equiv 0$ on
  $\R$. Hence, $f(x)\equiv 0$ on $\R$, which contradicts the
  condition. The first inequality is proved. The second inequality
  follows from the first one, if we take
  $g(x):=e^{-\alpha|x|}(1-\beta|x|)$. Then, for the indicated values
  of the parameters, $g\in C(\R)\cap L(\R)$ and $\widehat{g}(t)=
  \frac{2((\alpha-\beta)\alpha^2+(\alpha+\beta)t^2)}
  {(\alpha^2+t^2)^2}\ge 0$ for all $t\in\R$.
\end{proof}

\begin{z_lemma}\label{le4}
  Let $\nu$ be a function of bounded variation on a segment
  $[0,\sigma]$. Then $
  \lim_{t\to+\infty}\int_{0}^{\sigma}e^{-tu}d\nu(u)= \nu(+0)-\nu(0)
  $.
\end{z_lemma}

\begin{proof}
  Let first the function $\nu$ be right continuous in the point
  $t=0$. Then for any $\varepsilon\in(0,\sigma)$ and $t>0$, we have
  $ \left| \int_{0}^{\sigma}e^{-tu}d\nu(u) \right|\le
  V_{0}^{\varepsilon}+e^{-\varepsilon t}V_{0}^{\sigma} $. Here
  $V_{0}^{t}$ is the variation of the function $\nu$ on the segment
  $[0,t]$. By passing to the limit, we get that $
  \limsup_{t\to+\infty}\left| \int_{0}^{\sigma}e^{-tu}d\nu(u)
  \right|\le V_{0}^{\varepsilon} $. Since $\nu$ is right continuous
  in the point $t=0$, we have $ \lim_{\varepsilon\to+0}
  V_{0}^{\varepsilon}=0 $. In this case, the lemma is proved. In the
  general case, we assume that $\nu_1(0):=\nu(+0)$ and
  $\nu_1(t):=\nu(t)$ for $0<t\le\sigma$. It is clear that $\nu_1$ is
  a function of bounded variation on the segment $[0,\sigma]$ and is
  right continuous in the point $t=0$. Then $
  \int_{0}^{\sigma}e^{-tu}d\nu(u)=
  \int_{0}^{\sigma}e^{-tu}d\nu_1(u)+\nu(+0)-\nu(0) \to
  \nu(+0)-\nu(0)$ for $t\to+\infty$.
\end{proof}

\subsection{Statements about Functions in the Class $\overline{HB}$}

The following properties were proved in~\cite[Ch. VII]{Lev1}.
\begin{itemize}{\it
  \item[{\bf 1)}] If $\omega(z)\in\overline{HB}$, then the common
    zeros, if they exist, of the functions $\omega(z)$ and
    $\overline{\omega}(z)$ are real.

  \item[{\bf 2)}] Trivial functions in the class $\overline{HB}$ do
    not have zeros in $\C\setminus\R$.

  \item[{\bf 3)}] Let a function $\omega(z)$ be nontrivial. Then
    $\omega(z)\in\overline{HB}\iff \omega(z)=R(z)\omega_1(z)$, where
    $R(z)$ is a real entire function that does not have zeros in
    $\C\setminus\R$, and $\omega_{1}(z)\in HB$.

  \item[{\bf 4)}] If a sequence of functions
    $\omega_n(z)\in\overline{HB}$ converges to a function
    $\omega(z)\not\equiv 0$ on every compact subset of $\C$, then
    $\omega(z)\in\overline{HB}$.  }
\end{itemize}

For a function $\omega(z)=P(z)+iQ(z)$, where $P(z)$ and $Q(z)$ are
real entire functions, set $d(z):=P(z)Q'(z)- P'(z)Q(z)$ and $
H_{\alpha}(z):=P(z)\times\break\cos{\alpha} - Q(z)\sin{\alpha}$.
If the function $\omega(z)$ is real, up to a constant factor, then
it is clear that $d(x)\equiv 0$.

\begin{oldtheorem}
  [{B.~Ya.~Levin \cite[Ch.~VII, Th.~4]{Lev1} and N.~N.~Mejman
    \cite[Ch. IV, Th.~15$'$]{Cheb_Me}}] \label{old1} $\omega(z)\in
  HB\iff$ {\bf 1)} the functions $P(z)$ and $Q(z)$ do not have
  common zeros and, for any $\mu,\nu\in\R$, $|\mu|+|\nu|\ne 0$, the
  function $\mu P(z)+\nu Q(z)$ does not have zeros in
  $\C\setminus\R$; {\bf 2)} for some $x_0\in\R$, we have $d(x_0)>0$.
  Moreover, if conditions {\rm 1)} and {\rm 2)} hold, then $d(x)>0$
  holds for any $x\in\R$.
\end{oldtheorem}
\begin{oldtheorem}[{B.~Ya.~Levin  \cite[Ch.~VII, Th.~4$'$]{Lev1}}]
  \label{old2}
  Let a function $\omega(z)$ be nontrivial. Then $\omega(z)\in
  \overline{HB}\iff$ {\bf 1)} for any $\mu,\nu\in\R$,
  $|\mu|+|\nu|\ne 0$, the function $\mu P(z)+\nu Q(z)$ does not have
  zeros in $\C\setminus\R$; {\bf 2)} for some $x_0\in\R$, we have
  $d(x_0)>0$. Moreover, if conditions {\rm 1)} and {\rm 2)} are
  satisfied, then the inequality $d(x)\ge 0$ holds for any $x\in\R$.
 \end{oldtheorem}
 The following proposition immediately follows from
 Theorems~\ref{old1} and \ref{old2}.
  \pagebreak
\begin{z_proposition}\label{pr1}
  Let a function $\omega(z)\in \overline{HB}$ be nontrivial. Then
  \begin{itemize}
  \item[{\bf 1)}] $d(x)\ge 0$, $x\in\R$.

  \item[{\bf 2)}] $d(x_0)= 0$ for some $x_0\in\R \iff \omega(x_0)=0
    \iff P(x_0)=Q(x_0)=0$. If a number $x_0\in\R$ is a zero of
    multiplicity $p$ for the function $\omega$, then $x_0$ is a zero
    of the function $d$ of multiplicity $2p$.

  \item[{\bf 3)}] For all $\alpha\in\R$, the function $H_{\alpha}$
    is real and does not have zeros in $\C\setminus\R$. If a number
    $x_0\in\R$ is a zero of the function $H_{\alpha}$ of
    multiplicity $q$, then $q\le p+1$, where $p$ is the multiplicity
    of the zero $x_0$ for the function $\omega$ $(p=0$, if
    $\omega(x_0)\ne 0)$. If the function $\omega$ does not have real
    zeros, then all zeros of the function $H_{\alpha}$, if it has
    any, are simple.
  \end{itemize}
\end{z_proposition}

\begin{proof}
  By property~3, $\omega(z)=R(z)\omega_1(z)$, where $R(z)$ is a real
  entire function that does not have zeros in $\C\setminus\R$, and
  $\omega_{1}(z)\in HB$. Then $\omega(x_0)= 0$ for some $x_0\in\R
  \iff R(x_0)=0$ and, in this case, multiplicities of the zero $x_0$
  for $\omega(z)$ and $R(z)$ coincide.  If
  $\omega_{1}(z)=P_{1}(z)+iQ_{1}(z)$, where $P_{1}(z)$ and
  $Q_{1}(z)$ are real entire functions, then, by Theorem~\ref{old1},
  $d_{1}(x):=P_{1}(x)Q'_{1}(x)- P'_{1}(x)Q_{1}(x)>0$, $x\in\R$. It
  is clear that $P(z)=R(z)P_{1}(z)$ and $Q(z)=R(z)Q_{1}(z)$. Then
  $d(x)=R^2(x)d_1(x)$. Hence, $d(x_0)= 0$ for some $x_0\in\R \iff
  R(x_0)=0$ and, in this case, the multiplicity of the zero $x_0$ of
  the function $d$ is two times greater than the multiplicity of the
  zero $x_0$ of the function $R$. Statements~1) and 2) are proved.

  To prove statement~3) it should be noted that, by
  Theorem~\ref{old1}, the function $
  H_{1,\alpha}(z):=P_{1}(z)\cos{\alpha} - Q_{1}(z)\sin{\alpha}$ does
  not have zeros in $\C\setminus\R$ for all $\alpha\in\R$, and all
  its real zeros, if there are any, are simple. This follows from
  the identity
  $H_{1,\alpha}(x)H'_{1,\beta}(x)-H'_{1,\alpha}(x)H_{1,\beta}(x)\equiv
  d_{1}(x)\sin(\alpha-\beta)$. It remains to make a use of
  $H_{\alpha}(x)=R(x)H_{1,\alpha}(x)$.
\end{proof}

An indicator of growth of an exponential type function is defined by
$
h_{f}(\varphi):=\limsup_{r\to+\infty}\frac{\ln|f(re^{i\varphi})|}{r}$,
$ \varphi\in\R $.

\begin{z_definition}\label{z:d3}
  A function $\omega(z)$ is called a class $P$ function if it is an
  exponential type function, does not have zeros in the open lower
  half-plane $\Im z<0$ and
  $2d_{\omega}:=h_{\omega}\left(-\frac{\pi}{2}\right)-h_{\omega}\left(\frac{\pi}{2}\right)\ge
  0$ $($the quantity $d_{\omega}$ is called a defect of the function
  $\omega)$.
\end{z_definition}

\begin{oldtheorem}[{B.~Ya.~Levin \cite[Ch.~VII, Lemma~1]{Lev1}}]
  \label{old3}
  $\omega(z)\in P\iff \omega(z)\in\overline{HB}$ and $\omega(z)$ is
  an exponential type entire function.
\end{oldtheorem}

Since $\overline{\omega}(z)\equiv\overline{\omega(\overline{z})}$,
it is clear that the product of two class $HB$ functions is also a
class $HB$ function, that is, $HB\cdot HB\subset HB$. Similarly,
$\overline{HB}\cdot \overline{HB}\subset \overline{HB}$. It
follows from Theorem~\ref{old3} that $P\cdot P\subset P$, too. The
function classes $HB$ and $P$ were introduced and studied,
correspondingly, by M.~G.~Krein and B.~Ya.~Levin. The given
Definition~\ref{z:d1} is due to N.~N.~Meiman

Let $\mu(t)$ be a function with bounded variation on the segment
$[a,b]$, $a<b$, which is left continuous in every point of the
interval $(a,b)$, and $ \omega(z):=\int_{a}^{b} e^{izt}d\mu(t)$.

The following results of such functions are contained in the
mono\-graph~\cite[Ch.~I]{Leont}; they allow to easily determine the
defect. Let $[a_1,b_1]$ be the smallest segment that is contained in
$[a,b]$ and possessing the following property: the function $\mu(t)$
is constant on $[a,a_1]$ and $(b_1,b]$. If there is no such
intervals $[a,a_1]$ or $(b_1,b]$, then we take $a_1=a$ or $b_1=b$,
correspondingly. If $a_1=b_1$, then $\omega(z)\equiv ce^{ia_1 z}$.
If $a_1<b_1$, then $\omega(z)=\int_{a_1}^{b_1} e^{izt}d\mu_1(t)$,
where the function $\mu_1(t)$ coincides with $\mu(t)$ for $a_1\le
t<b_1$ and $\mu_1(b_1):=\mu(b)$. In this case, see~\cite[Ch. I, \S
4.3]{Leont}, the function $\omega(z)$ has infinitely many zeros and
$h_{\omega}\left(-\frac{\pi}{2}\right)=b_1$,
$h_{\omega}\left(\frac{\pi}{2}\right)=-a_1$ and, hence,
$2d_{\omega}=b_1+a_1$. Moreover, the upper limit in the definition
of the growth indicator, for almost all $\varphi\in\R$, is equal to
the limit.

\begin{z_proposition}\label{pr2}
  Let $F(z):=\int_{0}^{\sigma} e^{izt}d\mu(t)$, where $\mu(t)$ is a
  function with bounded variation on the segment $[0,\sigma]$,
  $\sigma>0$. Then,
  \begin{itemize}
  \item[{\bf 1)}] if $F(z)\not\equiv ce^{i\alpha z}$,
    $\alpha\in[0,\sigma]$, then the function $F$ has infinitely many
    zeros;
  \item[{\bf 2)}] $F\in\overline{HB}\iff F$ does not have zeros in
    the open lower half-plane $\Im z<0$.
  \end{itemize}
\end{z_proposition}

\begin{proof}
  Without loss of generality, we can assume that the function
  $\mu(t)$ is left continuous in every point of the interval
  $(0,\sigma)$. Statement~1) has been considered above. Necessity
  of~2) is immediate. Let us prove sufficiency of~2). Assume that
  the function $F$ does not have zeros in the open lower half-plane
  $\Im z<0$. If $F(z)\equiv ce^{i\alpha z}$, then $c\ne 0$ and
  $\alpha\in[0,\sigma]$, by Lemma~\ref{le1}, statement~4). In this
  case, it is easy to check that $F\in\overline{HB}$, and if
  $\alpha>0$, then $F\in HB$. If $F(z)\not\equiv ce^{i\alpha z}$,
  then, by the above, $0\le a_1<b_1\le \sigma$,
  $2d_{\omega}=b_1+a_1>0$ and, hence, $F\in P$. By
  Theorem~\ref{old3}, $F\in \overline{HB}$.
\end{proof}

\fussy
\subsection{An Interpolation Formula}

Denote by $B_{\sigma}^m$, $m\in\Z_+:=\N\cup\{0\}$, the class of
functions $f\in E_{\sigma}$, for which $f(x)=o(x^m)$,
$x\to\pm\infty$, on the real axis.

By $S_{\sigma}$, $\sigma>0$, denote the class of sine-type
functions, that is, the set of functions $F\in E_{\sigma}$ that
satisfy the following conditions: {\it {\bf
    a)}~$h_F(\pm\frac{\pi}{2})=\sigma$; {\bf b)} all zeros
  $\lambda_k$ of the function $F$ are simple and satisfy the
  condition $\inf_{k\ne n}|\lambda_k-\lambda_n|=2\delta>0$; {\bf c)}
  all zeros are located in a strip parallel to the real axis, that
  is, $\sup_{k}|\Im \lambda_k |=H<\infty$; {\bf d)}~there exist
  constants $C_k ,h\in\R$, such that $0<C_1\le |F(x+ih)|\le
  C_2<\infty$, $x\in\R$.}

If $F\in S_{\sigma}$, then $\sigma (F)=\sigma>0$ and $F$ has
infinitely many zeros both in the left half-plane $\Re z\le 0$ and
in the right half-plane $\Re z\ge 0$. Zeros of the function $F\in
S_{\sigma}$ are aways indexed in the increasing order of their real
parts, that is, $\Re \lambda_{k}\le \Re \lambda_{k+1}$, $k\in\Z$.
The function $F(z):=\sin(\sigma z+\alpha)$ is an example of such a
function.

\begin{oldtheorem}[\!\!\!{\cite[Lemma~1]{Zast2004_MFAT}}]
  \label{old5}
  Let $F\in S_{\sigma}$, $\sigma >0$, and $\{\lambda_k\}$ be a
  sequence of all its zeros. Then for any $m\in\Z_+$, $f\in
  B_{\sigma}^m$, $\tau\in\C$, and $z\in\C ,z\ne \lambda_k+\tau $, we
  have
  \begin{equation*} \label{old}
    \frac{d^m}{du^m} \Bigl\{\frac{f(u)}{F(u-\tau)}
    \Bigr\}\Big|_{u=z}
    = - m! \lim_{n\to\infty}\sum_{|\lambda_k|<n} \frac{f(\lambda_k+\tau)}{F'(\lambda_k)
      (\lambda_k+\tau-z)^{m+1}}.
  \end{equation*}
\end{oldtheorem}

Note that, for a smaller class of functions $f$, this formula is well-known, see details in~\cite[\S~1]{Zast2004_MFAT}.

\section{Proof of Theorems~\ref{z:t1} and~\ref{z:t2}}

\begin{proof}[Proof of Theorem~\ref{z:t1}]
  If $f$ is an entire function of exponential type $\le\sigma$,
  $\sigma>0$, and $f(x)=o(x)$, $x\to\pm\infty$, then the following
  interpolation formula holds for any $\alpha$ and $x$:
  \begin{multline}\label{b1}
    \sigma f(x)\cos(\sigma x+\alpha)-f'(x)\sin(\sigma x+\alpha)\\=
    \sigma\lim_{n\to+\infty}\sum_{k=-n}^{n} \frac{\sin^2(\sigma
      x+\alpha)}{(\sigma x+\alpha-k\pi)^2}\cdot(-1)^k
    f\Bigl(\frac{k\pi-\alpha}{\sigma}\Bigr).
  \end{multline}
  This follows from Theorem~\ref{old5} for $F(z):=\sin(\sigma
  z+\alpha)$, $\lambda_k=\frac{k\pi-\alpha}{\sigma}$, $m=1$,
  $\tau=0$, $z=x$.

  Apply formula~(\ref{b1}) to the function
  $f(x):=P\left(x-\frac{\tau}{\sigma}\right)
  \cos\alpha-\break Q\left(x-\frac{\tau}{\sigma}\right)\sin\alpha$,
  $\alpha\in\R$. Since
  $(-1)^kf\left(\frac{k\pi-\alpha}{\sigma}\right)=
  E\left(\frac{k\pi-\alpha-\tau}{\sigma}\right)\ge 0$ for all
  $k\in\Z$, we have
  \begin{multline}\label{bb2}
    \sigma\left(P\Bigl(x-\frac{\tau}{\sigma}\Bigr)\cos\alpha-Q\Bigl(x-\frac{\tau}{\sigma}\Bigr)\sin\alpha\right)\cos(\sigma
    x+\alpha)
    \\
    -\left(P'\Bigl(x-\frac{\tau}{\sigma}\Bigr)\cos\alpha-Q'\Bigl(x-\frac{\tau}{\sigma}\Bigr)\sin\alpha\right)\sin(\sigma
    x+\alpha) =
    \\
    \sigma\sum_{k=-\infty}^{+\infty} \frac{\sin^2(\sigma
      x+\alpha)}{(\sigma x+\alpha-k\pi)^2}\cdot
    E\Bigl(\frac{k\pi-\alpha-\tau}{\sigma}\Bigr) \ge 0,\qquad
    \alpha,x\in\R.
  \end{multline}
  Consider first the case $\tau=0$. Then
  \begin{multline}\label{b2}
    \sigma(P(x)\cos\alpha-Q(x)\sin\alpha)\cos(\sigma x+\alpha)\\-
    (P'(x)\cos\alpha-Q'(x)\sin\alpha)\sin(\sigma x+\alpha)\ge
    0,\qquad\alpha,x\in\R.
  \end{multline}
  Assume that the identities
  $E\bigl(\frac{k\pi-\beta}{\sigma}\bigr)= 0$ hold for some
  $\beta\in\R$ and all $k\in\Z$. Then inequality~(\ref{b2}) becomes
  identity in $x\in\R$ for $\alpha=\beta$ and, hence,
  $P(x)\cos\beta-Q(x)\sin\beta\equiv \gamma\sin(\sigma x+\beta)$,
  $x\in\R$, for some constant $\gamma\in\R$. Let, for example,
  $\cos\beta\ne 0$. Then expressing $P$ in terms of $Q$ and
  substituting it into~(\ref{6}) for $E$ we get the identity $
  E(x)\cos\beta\equiv f_1(x)\sin(\sigma x+\beta)$, where
  $f_1(x):=\gamma\cos\sigma x+Q(x)$. Since $E(x)\ge 0$ for all
  $x\in\R$, all real zeros of the function $E$ has even
  multiplicity. Hence, $f_1\bigl(\frac{k\pi-\beta}{\sigma}\bigr)= 0$
  for all $k\in\Z$. Applying formula~(\ref{b1}) to the function
  $f_1$, we get the identity $ \gamma\cos\sigma x+Q(x)\equiv
  c_1\sin(\sigma x+\beta)$ for some constant $c_1\in\R$. Setting
  $c:=\frac{c_1}{\cos\beta}$ we get identities~(\ref{8}). In a
  similar way, we can consider the case $\sin\beta\ne 0$. One can
  directly check that inequality~(\ref{7}) becomes identity
  if~(\ref{8}) holds and, in this case, we get
  $d(x)\equiv\gamma^2\sigma$.

  Assume now that, for all $\alpha\in\R$, $E(x)\not\equiv
  c\sin^2(\sigma x+\alpha)$. Then for any $\alpha\in\R$ there exists
  $k_0\in\Z$ such that $ E\bigl(\frac{k_0\pi-\alpha}{\sigma}\bigr)>
  0$.  In this case, inequality~(\ref{b2}) is strict for all
  $x\ne\frac{k\pi-\alpha}{\sigma}$, $k\in\Z$. Hence,
  \begin{equation}\label{b3}
    \text{
      \parbox{.85\textwidth}
      {
        inequality~(\ref{b2}) becomes identity
        for some $x=x_0\in\R$ and $\alpha=\alpha_0\in\R$
        $\iff$ $x_0=\frac{k_0\pi-\alpha_0}{\sigma}$ and
        $E(x_0)=0$ for some $k_0\in\Z$.
      }
    }
  \end{equation}
  Let
  \begin{equation*}
    \begin{aligned}
      A_1(x)&:= \sigma P(x)\cos\sigma x-P'(x)\sin\sigma x ,\\
      A_2(x)&:= \sigma Q(x)\sin\sigma x+Q'(x)\cos\sigma x , \\
      A_3(x)&:=\sigma (P(x)\sin\sigma x+Q(x)\cos\sigma x)+
      P'(x)\cos\sigma x-Q'(x)\sin\sigma x .
    \end{aligned}
  \end{equation*}
  Then inequality~(\ref{b2}) is equivalent to the inequality
  \begin{multline}\label{b4}
    A_1(x)+A_2(x)+(A_1(x)-A_2(x))\cos 2\alpha\\-A_3(x)\sin 2\alpha\ge
    0,\qquad\alpha,x\in\R.
  \end{multline}
  Inequality~(\ref{b4}) with two parameters is equivalent to the
  following inequality with one parameter:
  \begin{equation}\label{b5}
    \sqrt{(A_1(x)-A_2(x))^2+A^2_3(x)}\le A_1(x)+A_2(x) ,\qquad x\in\R.
  \end{equation}
  Here, see~(\ref{b3}),
  \begin{equation*}
    \text{
      \parbox{.91\textwidth}
      {
        inequality~(\ref{b5}) becomes identity for some
        $x=x_0\in\R$ $\iff$
        inequality~(\ref{b4}) becomes identity for
        $x=x_0\in\R$ and some $\alpha=\alpha_0\in\R\iff E(x_0)=0$.
      }
    }
  \end{equation*}
  Since $A_1(x)\ge 0$ for all $x\in\R$ (this is
  inequality~(\ref{b2}) with $\alpha=0$) and $A_2(x)\ge 0$ for all
  $x\in\R$ (this is inequality~(\ref{b2}) with
  $\alpha=\frac{\pi}{2}$), inequality~(\ref{b5}) is equivalent to
  the inequality
  \begin{equation}\label{b6}
    A^2_3(x)\le 4A_1(x)A_2(x) ,\qquad x\in\R\;.
  \end{equation}
  Inequality~(\ref{b6}) is equivalent to inequality~(\ref{7}). This
  is implied by the following identity:
  \begin{multline*}
    A^2_3(x)- \left\{(\sigma P(x)+Q'(x))\sin\sigma x+(P'(x)-\sigma
      Q(x))\cos\sigma  x\right\}^2 \\\equiv
    4A_1(x)A_2(x)-4\sigma\left(P(x)Q'(x)-P'(x)Q(x)\right) \;.
  \end{multline*}
  Hence, if $\tau=0$, statements~1), 2)~$iii)\Rightarrow
  iv\Rightarrow i)$ and~3) are proved. The general case is reduced
  to the previous one by considering the functions
  $P_1(x):=P\left(x-\frac{\tau}{\sigma} \right)$ and
  $Q_1(x):=Q\left(x-\frac{\tau}{\sigma} \right)$. Then
  $E_1(x):=P_1(x)\cos\sigma x+Q_1(x)\sin\sigma x=
  E\left(x-\frac{\tau}{\sigma} \right)\ge 0$, $x\in\R$.

  Let us prove the remaining statements in~2). The implication
  $ii)\Rightarrow iii)$ is clear. Let inequality~(\ref{7}) become an
  identity. Assume that, for any $\beta\in\R$ and $c\ge 0$,
  $E(x)\not\equiv c\sin^2(\sigma x+\tau+\beta)$ and, hence,
  $E(x)\not\equiv 0$. But it follows from statement~3) that
  $E(x)\equiv 0$. This contradiction proves the implication
  $i)\Rightarrow ii)$. Theorem~\ref{z:t1} is proved.
\end{proof}

\begin{z_remark}\label{re1}
  It can be seen from the above proof that the following converse
  proposition also holds. {\it Let $\omega(z)=P(z)+iQ(z)$, where
    $P,Q$ are real functions of the class $E_{\sigma}$, $\sigma >
    0$, and $\omega(x)=o(x)$, $x\to\pm\infty$, on the real axis. If
    for some $\tau\in\R$ and all $x\in\R$ the identity~$(\ref{7})$
    holds and, moreover, the inequalities $A_1(x):= \sigma
    P(x)\cos(\sigma x+\tau)-P'(x)\sin(\sigma x+\tau)\ge 0$ and
    $A_2(x):= \sigma Q(x)\sin(\sigma x+\tau)+Q'(x)\cos(\sigma
    x+\tau) \ge 0$ hold for all $x\in\R$, then the inequality $
    E(x):= P(x)\cos(\sigma x+\tau)+Q(x)\sin(\sigma x+\tau)\ge 0$,
    $x\in\R$ also holds.  }
\end{z_remark}

\begin{z_proposition}\label{pr3}
  Let the conditions of Theorem~{\rm \ref{z:t1}} be satisfied and
  $H_{\alpha}(z):=P(z)\cos{\alpha} - Q(z)\sin{\alpha}$. Then the
  following conditions are equivalent:
  \begin{itemize}
  \item[{\bf i)}] the function $\omega$ is real up to a constant
    factor;
  \item[{\bf ii)}] $d(x)\equiv 0$;
  \item[{\bf iii)}] for some $c\ge 0$, $\beta\in\R$, the identity
    $\omega(x)\equiv ce^{i(\frac{\pi}{2}-\beta)}\sin(\sigma
    x+\tau+\beta)$ holds;
  \item[{\bf iv)}] for some $\alpha\in\R$, we have
    $H_\alpha(x)\equiv 0$.
  \end{itemize}
\end{z_proposition}

\begin{proof}
  Let us prove the implication~i) $\Rightarrow$ ii). Let
  $\omega(z)\equiv e^{i\beta}\omega_0(z)$, where $\omega_0$ is
  real and $\beta\in\R$. Then $P(x)=\omega_0(x)\cos\beta$,
  $Q(x)=\omega_0(x)\sin\beta$ and, clearly, $d(x)\equiv 0$.

  ii) $\Rightarrow$ iii). Let $d(x)\equiv 0$. Then
  inequality~(\ref{7}) becomes identity. It follows from
  Theorem~\ref{z:t1} that, for some $c\ge 0$ and $\gamma\in\R$,
  identities~(\ref{8}) hold and, moreover, $d(x)\equiv
  \gamma^2\sigma$. Hence, $\gamma=0$ and so $\omega(x)\equiv
  ce^{i(\frac{\pi}{2}-\beta)}\sin(\sigma x+\tau+\beta)$.

  iii) $\Rightarrow$ iv). Let the identity $\omega(x)\equiv
  ce^{i(\frac{\pi}{2}-\beta)}\sin(\sigma x+\tau+\beta)$ hold for
  some $c\ge 0$, $\beta\in\R$. Then $P(x)=c\sin\beta \sin(\sigma
  x+\tau+\beta)$, $Q(x)=c\cos\beta \sin(\sigma x+\tau+\beta)$ and,
  hence, $H_{\alpha}(x)=c\sin(\sigma x+\tau+\beta)(\sin\beta
  \cos\alpha -\cos\beta \sin\alpha )\equiv 0$ for $\alpha=\beta$.

  iv) $\Rightarrow$ i). Let, for some $\alpha\in\R$,
  $H_{\alpha}(x)\equiv 0$. Then $ P(x)\cos{\alpha} -
  Q(x)\sin{\alpha}\equiv 0$. Hence, either $Q(x)\equiv\lambda P(x)$,
  or $P(x)\equiv\lambda Q(x)$ for some $\lambda\in\R$. In any case,
  $\omega$ is a real function up to a constant factor.
  Proposition~\ref{pr3} is proved.
\end{proof}

\begin{z_proposition}\label{pr4}
  Let the conditions of Theorem~{\rm \ref{z:t1}} be satisfied,
  $H_{\alpha}(z)\break:=P(z)\cos{\alpha} - Q(z)\sin{\alpha}$, and assume
  that the function $\omega$ is not real up to a constant multiple.
  Then we have the following.
  \begin{itemize}
  \item[{\bf i)}] $H_{\alpha}$, for any $\alpha\in\R$, is a real
    function of the class $E_{\sigma}$, $H_{\alpha}\not\equiv 0$,
    $H_{\alpha}(x)=o(x)$, $x\to\pm\infty$, on the real axis, and
    $(-1)^pH_{\alpha}\left( \frac{p\pi-\alpha-\tau}{\sigma}\right)\break=
    E\left(\frac{p\pi-\alpha-\tau}{\sigma}\right)\ge 0$, $p\in\Z$.
  \item[{\bf ii)}] For any $\alpha\in\R$, the function $H_{\alpha}$
    has infinitely many zeros and all of them are real,
    $xH_{\alpha}(x)\not=o(1)$, $x\to\pm\infty$, on the real axis. In
    every interval $I_p:=(\lambda_{p-1},\lambda_{p})$, where
    $\lambda_p=\lambda_p(\alpha):=\frac{p\pi-\alpha-\tau}{\sigma}$,
    the function $H_{\alpha}$ can have only one zero, and it is
    simple. Moreover, if $x_0\in I_p$ and $H_{\alpha}(x_0)=0$, then
    $(-1)^pH'_{\alpha}(x_0)>0$. If, for some $p\in\Z$, the number
    $\lambda_{p}$ is a zero of the function $H_{\alpha}$, then its
    multiplicity does not exceed $2$, and one of the intervals $I_p$
    or $I_{p+1}$ does not contain zeros of the function
    $H_{\alpha}$. If the number $\lambda_{p}$ is a zero of the
    function $H_{\alpha}$ of multiplicity $2$, then
    $(-1)^pH^{(2)}_{\alpha}(\lambda_{p})<0$ and
    $(-1)^pH_{\alpha}(x)<0$ for $x\in I_p\cup I_{p+1}$, and the
    numbers $\lambda_{p-1}$ and $\lambda_{p+1}$ could only be simple
    zeros.
  \item[{\bf iii)}] We have $\omega\in\overline{HB}$.
  \item[{\bf iv)}] If the function $\omega$ has real zeros, then
    they are simple.
  \item[{\bf v)}] $d(x_0)= 0$ for some $x_0\in\R \iff \omega(x_0)=0
    \iff P(x_0)=Q(x_0)=0 $. If a number $x_0\in\R$ is a zero of the
    function $\omega$, then the number $x_0$ is a zero of the
    function $d$ of multiplicity $2$.
  \item[{\bf vi)}] If a number $x_0\in\R$ is a multiplicity $2$ zero
    of the function $H_{\alpha}$, then $\omega(x_0)=0$.
  \item[{\bf vii)}] If $E(x)>0$, $x\in\R$, then $\omega\in HB$, and
    all zeros of the function $H_{\alpha}$, $\alpha\in\R$, are
    simple.
\end{itemize}
\end{z_proposition}

\begin{proof}
  Assume that the function $\omega$ is not real up to a constant
  factor. Then $H_{\alpha}(x)\not\equiv 0$ for any $\alpha\in\R$
  (Proposition~\ref{pr3}). The remaining part of statement~i) is
  clear. Statement~ii) follows from statement~i) and
  Theorem~\ref{old4}.

\begin{oldtheorem}[\!\!\!{\cite[Theorem~1 for $F(z):=\sin(\sigma z+\beta),\ \beta\in\R,$ $\lambda_k=\frac{k\pi-\beta}\sigma$,
    $n=0$]{Zast2004_MFAT}\footnote{For a smaller class of functions
      $f$, this theorem was proved in author's
      work~\cite{Zast2004_MZ}.}}]
  \label{old4}
  Let a function $f$ satisfy the following conditions: {\bf a)} $f $
  is a real entire function of exponential type $\le\sigma$,
  $\sigma>0$, $f\not \equiv 0$, and $f(x)=o(x)$, $x\to\pm\infty $ on
  the real axis; {\bf b)}~for some $\beta\in \R$ and $\forall k\in
  \Z$, we have $(-1)^{k}f(\lambda_k)\ge 0$, where
  $\lambda_k:=\frac{k\pi-\beta}{\sigma}$. Then the following holds.
  \begin{itemize}
  \item[{\bf 1)}] In every interval
    $I_p:=(\lambda_{p-1},\lambda_{p})$, $p\in \Z$, there may exist
    only one zero of the function $f$ and if there is one, then it
    is simple. Moreover, if $x_0\in I_p$ and $f(x_0)=0$, then
    $(-1)^pf'(x_0)>0$.\footnote{This inequality is contained in the
      proof of this theorem.}
  \item[{\bf 2)}] $xf(x)\ne o(1),x\to\pm\infty.$
  \item[{\bf 3)}] The function $f$ has only real zeros.
  \item[{\bf 4)}] If for some $p\in \Z$, the number $\lambda_{p}$ is
    a zero of the function $f,$ then its multiplicity is not greater
    than $2$, and one of the intervals $I_p$ or $I_{p+1}$ does not
    contain zeros of the function $f$. If the number $\lambda_{p}$
    is a zero of multiplicity $2,$ then $(-1)^{p}
    f^{(2)}(\lambda_{p})<0$ and $(-1)^{p} f(x)<0$ for $x\in I_p \cup
    I_{p+1},$ and the numbers $\lambda_{p-1}$ and $\lambda_{p+1}$
    can only be simple zeros.
\end{itemize}
\end{oldtheorem}
Let us prove statement~iii). It follows from inequality~(\ref{7})
and Proposition~\ref{pr3} that $d(x_0)>0$ for some $x_0\in\R$.
Statement~ii) implies that, for any $\alpha\in\R$, the function
$H_{\alpha}$ does not have zeros in $\C\setminus\R$. Since, by the
condition, the function $\omega$ is not trivial, by
Theorem~\ref{old2}, $\omega\in\overline{HB}$.

Let us prove statement~iv). Assume that for some $x_0\in\R$, we have
$\omega(x_0)=\omega'(x_0)=0$. Then
$P(x_0)=Q(x_0)=P'(x_0)=Q'(x_0)=0$. Thus, for any $\alpha\in\R$, the
number $x=x_0$ is a zero of the function $H_{\alpha}$, and its
multiplicity is not less that $2$. Statement~ii) implies that
$x_0\in\bigl\{\frac{p\pi-\alpha-\tau}{\sigma}:p\in\Z\bigr\}\cap\bigl\{\frac{p\pi-\delta-\tau}{\sigma}:p\in\Z\bigr\}$,
$\alpha,\delta\in\R$, however, for $\alpha=\delta-\frac{\pi}{2}$,
this intersection is empty.

Let us prove statements~v) and~vi). As was proved in~iii),
$\omega\in\overline{HB}$. But then, we can apply
Proposition~\ref{pr1} and use statement~iv).

Now we prove~vii). It is clear that all real zeros of the function
$\omega$ are zeros of the function $E$. Hence, if $E(x)>0$,
$x\in\R$, then the function $\omega$ does not have real zeros, and
all zeros of the function $H_{\alpha}$, $\alpha\in\R$, are real, see
statement~ii), and simple, see Proposition~\ref{pr4}.
\end{proof}

\begin{z_proposition}\label{pr5}
  Assume that conditions of Theorem~{\rm \ref{z:t1}} hold. Then, for
  any $n\in\N$, the function
  $\omega^{(n)}(z)=P^{(n)}(z)+i\,Q^{(n)}(z)$ satisfies conditions of
  Theorem~{\rm \ref{z:t1}} for $\tau_n=\tau+\frac{\pi n}{2}$, that
  is, $P^{(n)},Q^{(n)}$ are real functions of class $E_{\sigma}$,
  $\omega^{(n)}(x)=o(x)$, $x\to\pm\infty$, on the real axis, and
  \begin{multline*} E_n(x):= P^{(n)}(x)\cos\left(\sigma
      x+\tau+\frac{\pi n}{2}\right)\\+Q^{(n)}(x)\sin\left(\sigma
      x+\tau+\frac{\pi n}{2}\right)\ge  0,\qquad x\in\R.
  \end{multline*}
  Moreover,
  \begin{itemize}
  \item[{\bf i)}] $E_n(x_0)=0$ for some $x_0\in\R\!\iff\!$
    $E(x)\equiv c\sin^2\bigl(\sigma x+\frac{\pi n}{2}-\sigma
    x_0\bigr)$ for some $c\ge 0$. In this case, inequality~\eqref{7}
    becomes identity for $\omega^{(n)}(z)$, and
    $d_n(x)\!:=\!P^{(n)}(x)Q^{(n+1)}(x)-
    P^{(n+1)}(x)Q^{(n)}(x)\equiv\gamma^2\sigma^{2n+1}$, where
    $\gamma$ is in~\eqref{8} with $\beta=\frac{\pi n}{2}-\tau-\sigma
    x_0$.
  \item[{\bf ii)}] $\omega^{(n)}$ is a real function up to a
    constant factor $\iff\omega$ is real up to a constant
    factor.
\end{itemize}
\end{z_proposition}

\begin{proof}
  Consider the case $n=1$. It was proved in~\cite[\S
  1]{Zast2004_MFAT} that $P',Q'\in E_{\sigma}$ and
  $\omega'(x)=o(x)$, $x\to\pm\infty$. If we take
  $\alpha=\frac{\pi}{2}-\sigma x$ in~(\ref{bb2}) and then replace
  $x$ with $x+\frac{\tau}{\sigma}$, we get the inequality
  $$
  E_1(x)= \sigma\sum_{k=-\infty}^{+\infty}
  \frac{E\left(\frac{k\pi-\frac{\pi}{2}+\sigma
        x}{\sigma}\right)}{(\frac{\pi}{2}-k\pi)^2} \ge 0,\qquad
  x\in\R.
  $$
  It follows at once from this inequality and Theorem~\ref{z:t1},
  statement~2) that~i) holds.

  Let us prove~ii). Let $\omega'$ be real up to a constant factor.
  It follows from Proposition~\ref{pr3} applied to $\omega'$ that,
  for some $c\ge 0$, $\beta\in\R$, we have the identity
  $\omega'(x)\equiv ce^{i(\frac{\pi}{2}-\beta)}\sin(\sigma
  x+\tau+\frac{\pi}{2}+\beta)$ and, hence, $\omega(x)\equiv
  c\sigma^{-1}e^{i(\frac{\pi}{2}-\beta)}\sin(\sigma
  x+\tau+\beta)+A+iB$, where $A,B\in\R$. Thus, $E(x)\equiv
  c\sigma^{-1}\sin^2(\sigma x+\tau+\beta)+A\cos(\sigma
  x+\tau)+B\sin(\sigma x+\tau)\ge 0$, $x\in\R$. Lemma~\ref{le1},
  statement~7 gives that $A=B=0$ and, hence, by
  Proposition~\ref{pr3} used for $\omega$, $\omega$ is real up to a
  constant. The converse is clear.

  Proposition~\ref{pr5} has been proved for $n=1$. The general case
  is proved by induction on $n\in\N$.
\end{proof}

\begin{proof}[Proof of Theorem~$\ref{z:t2}$]
  Statements~1) and~2) are proved in Proposition $\ref{pr4}$.

  Let us prove~3). Assume that a function $\omega$ is not real up to
  a constant. It follows from Propositions~\ref{pr5} and~\ref{pr4}
  that $\omega^{(n)}\in\overline{HB}$. Let us show that the function
  $\omega^{(n)}$ does not have real zeros. If for any $\alpha\in\R$
  and $c\ge 0$, $E(x)\not\equiv c\sin^2(\sigma x+\tau+\alpha)$, we
  have (see Proposition~\ref{pr5}) that $E_n(x)>0$, $x\in\R$, and,
  hence (see Proposition~\ref{pr4}.vii applied to $\omega^{(n)}$),
  the function $\omega^{(n)}$ does not have real zeros.

  Let, for some $c\ge 0$, $\beta\in\R$, the identity $E(x)\equiv
  c\sin^2(\sigma x+\tau+\beta)$ hold. Then, for some $\gamma\in\R$
  we have~(\ref{8}). Assume that $\omega^{(n)}(x_0)=0$ for some
  $x_0\in\R$. Then
  \begin{multline*}
    P^{(n)}(x_0)=c\sigma^n\sin\beta\sin\left(\sigma
      x_0+\tau+\beta+\frac{\pi
        n}{2}\right)\\[-2pt]+\gamma\sigma^n\sin\left(\sigma x_0+\tau+\frac{\pi
        n}{2}\right)=0,
  \end{multline*}\vskip-7mm
  \begin{multline*}
    Q^{(n)}(x_0)=c\sigma^n\cos\beta\sin\left(\sigma
      x_0+\tau+\beta+\frac{\pi
        n}{2}\right)\\[-2pt]-\gamma\sigma^n\cos\left(\sigma x_0+\tau+\frac{\pi
        n}{2}\right)=0.
  \end{multline*}
  Hence,
  \begin{multline*}
    P^{(n)}(x_0)\cos(\sigma x_0+\tau+\frac{\pi
      n}{2})+Q^{(n)}(x_0)\sin(\sigma x_0+\tau+\frac{\pi n}{2})\\
    = c\sigma^n\sin^2(\sigma x_0+\tau+\beta+\frac{\pi n}{2})=0,
  \end{multline*}
  and so $\gamma=0$. Then $\omega$ is real up to a constant factor,
  which contradicts the condition. Statement~3) is proved.
\end{proof}

\begin{z_remark}\label{re2}
  If conditions of Theorem~\ref{z:t1} are fulfilled and
  $\omega(z)\not\equiv 0$, it follows from Propositions~\ref{pr3}
  and~\ref{pr4}(ii) that $x\omega(x)\not=o(1)$, $x\to\pm\infty$, on
  the real axis. Moreover, Theorems~\ref{z:t1} and~\ref{z:t2} cease
  to hold if the condition $\omega(x)=o(x)$, $x\to\pm\infty$, is
  replaced with the condition $\omega(x)=O(x)$, $x\to\pm\infty$.
  This is easily seen by considering the function $\omega(z):=\sin
  z+az\cos z+i(az\sin z+1-\cos z)$, where $-1<a<-\frac{1}{2}$. It is
  clear that $\omega(x)=O(x)$ for $x\to\pm\infty$, and
  inequality~\eqref{6} holds for $\tau=-\frac{\pi}{2}$, $\sigma=1$;
  in this case, $E(x)=1-\cos x$. It is easy to check that
  $d(x)=a^2x^2+ax\sin x+(a+1)(1-\cos x)$. Since $d(x)\sim a^2x^2$
  for $x\to\pm\infty$ and $d(x)\sim x^2(a+1)(a+\frac{1}{2})$ for
  $x\to 0$, the function $d$ does not preserve the sign on the real
  axis and, hence, $\omega(z)\not\in \overline{HB}$; in the opposite
  case, Theorem~\ref{old2} implies that $d(x)\ge 0$, $x\in\R$. In
  the case Theorem~\ref{z:t2} is applied, we can consider the
  function $\omega(z):=zF(z)$, where the function $F$ is of the
  form~\eqref{1} and satisfies the conditions of Theorem~\ref{z:t4}
  (statement~2). Clearly, $\omega(x)=O(x)$ for $x\to\pm\infty$, and
  inequality~\eqref{6} holds for $\tau=-\frac{\pi}{2}$ (in this
  case, $E(x)=xS(x)\ge 0$, $x\in\R$), but $\omega(z)\not\in
  \overline{HB}$, since the function $F$ has one zero in the open
  lower half-plane.
\end{z_remark}

\section{Proofs of Theorems~\ref{z:t3} and~\ref{z:t4}}

In this section, we consider the functions defined by~\eqref{1},
\eqref{2}, \eqref{3}.

\begin{z_corollary}\label{sl2}
  Let $\mu$ be a real function with bounded variation on the segment
  $[0,\sigma]$, and let one of the four conditions hold,~\eqref{9}
  \eqref{10}, \eqref{11}, or~\eqref{12},
  \begin{equation}\label{9}
    C(x)\ge 0,\ x\in\R \ \text{ and }\ n=0,\ \tau=0,
  \end{equation}
  \begin{equation}\label{10}
    S(x)\ge 0 ,\ x>0,\  F(x)=o(1),\ x\to\pm\infty  \ \text{ and }\ n=1,\
    \tau=-\frac{\pi}{2},
  \end{equation}
  \begin{equation}\label{11}
    S(x)\ge 0  ,\ x>0,\   F(0)=0  \ \text{ and }\ n=-1,\
    \tau=-\frac{\pi}{2},
  \end{equation}
  \begin{equation}\label{12}
    \exists \,\tau_0\in\R:C(x)\cos{\tau_0}-S(x)\sin{\tau_0}\ge 0,\ x\in\R  \ \text{ and }\ n=0,\
    \tau=\pm\tau_0.
  \end{equation}
  Let $\omega(z):=z^n F(z)\equiv P(z)+iQ(z)$, where $P(z)\equiv z^n
  G(z)$, $Q(z)\equiv z^n H(z)$.  Then we have the following.
  \begin{itemize}
  \item[{\bf 1)}] The function $\omega$ satisfies the conditions of
    Theorem~$\ref{z:t1}$, that is, $1)$ $P,Q$ are real functions
    of class $E_{\sigma}$ and $\omega(x)=o(x)$, $x\to\pm\infty$, on
    the real axis; $2)$ for a corresponding value of $\tau$,
    $E(x):= P(x)\cos(\sigma x+\tau)+Q(x)\sin(\sigma
    x+\tau)\equiv x^n (C(x)\cos{\tau}-S(x)\sin{\tau})\ge 0,\
    x\in\R $,
    and, hence, for all $x\in\R$,
    \begin{equation}\label{13}
      4\sigma d(x)\ge x^{2n-2}D(x)\,,
    \end{equation}
    where
    \begin{multline*}
      D(x):=\bigl\{(2\sigma xS(x)+xC'(x)+nC(x))\cos\tau\\+(2\sigma
      xC(x)-xS'(x)-nS(x))\sin\tau\bigr\}^2
    \end{multline*}
    and $d(x):=P(x)Q'(x)-P'(x)Q(x)\equiv x^{2n}\Delta(x).$
  \item[{\bf 2)}] Inequality~\eqref{13} becomes identity $\iff
    E(x)\equiv 0$ $ \iff\!\! C(x)\cos\tau\equiv S(x)\sin\tau\equiv
    0$.  In this case, $F(z)\equiv\Delta(z)\equiv 0$ if $\cos\tau\ne
    0$, and $F(z)\equiv ce^{i\sigma z}$, $\Delta(z)\equiv
    c^2\sigma$, $c\in\R$ if $\cos\tau =0$.
  \item[{\bf 3)}] Inequality~\eqref{13} becomes identity for some
    $x=x_0\in\R\iff E(x_0)=0$.
  \item[{\bf 4)}] If $F(z)\not\equiv 0$, then the function $\omega$
    is not real up to a constant multiple and, hence,
    Theorem~$\ref{z:t2}$ and Proposition~$\ref{pr4}$ hold for the
    functions $\omega, P, Q, E, d$, and
    $H_{\alpha}(z):=P(z)\cos{\alpha} - Q(z)\sin{\alpha}\equiv z^n
    h_{\alpha}(z)$.
  \item[{\bf 5)}] If condition~\eqref{9} is satisfied and, in
    addition, $F(z)\not\equiv 0$, $S(x)\ge 0$ for $x>0$, $F(0)>0$,
    then the function $F$ does not have real zeros.
  \item[{\bf 6)}] If condition~\eqref{10} is satisfied and,
    additionally, $F(z)\not\equiv 0$, then $F(0)>0$, $H'(0)>0$, and
    $\Delta(0)>0$.
  \item[{\bf 7)}] If condition~\eqref{12} is satisfied and,
    additionally, $F(z)\not\equiv 0$, $\sin\tau_0\ne 0$, then the
    function $F$ does not have real zeros.
  \end{itemize}
\end{z_corollary}
\begin{proof}
  Statement~1) follows at once from Theorem~\ref{z:t1}, statement~1,
  if identities~(\ref{a2}) and~(\ref{a3}) are used.
  The first part in statement~2) immediately follows from
  Theorem~\ref{z:t1}, statement~2, and Lemma~\ref{le1}, statement~3.
  The second part of this statement follows from the first one and
  Lemma~\ref{le1}, statement~1.
  Statement~3) follows from statement~2) and Theorem~\ref{z:t1},
  statement~3.
  Statement~4) immediately follows from Lemma~\ref{le1}, statement~6.

  Let us now prove statement~5). If $x_0\in\R$ and $F(x_0)=0$, then
  $x_0\ne 0$ and $C(x_0)=S(x_0)=0$. Hence, $C'(x_0)=S'(x_0)=0$ and
  so, $F'(x_0)=0$, which contradicts statement~iv) in
  Proposition~\ref{pr4} about simplicity of real roots of the
  function $\omega(z)\equiv F(z)$. Statement~5) is proved.

  Now we prove statement~6). Let us first apply statement~ii) of
  Proposition~\ref{pr4} for $\alpha=0$. Then
  $\lambda_p=\frac{p\pi+\frac{\pi}{2}}{\sigma}$. The number
  $x_0=0\in I_0=(-\frac{\pi}{2},\frac{\pi}{2})$ and it is a zero of
  the function $H_0(x)=xG(x)$. Hence, $H'_0(0)>0$. It remains to use
  that $H'_0(0)=G(0)=F(0)$.
  Now apply statement~ii) of Proposition~\ref{pr4} for
  $\alpha=-\frac{\pi}{2}$. Then $\lambda_p=\frac{(p+1)\pi}{\sigma}$.
  The number $\lambda_{-1}=0$ is a zero of the function
  $H_{-\frac{\pi}{2}}(x)=xH(x)$ of multiplicity not less than $2$.
  Thus, $H''_{-\frac{\pi}{2}}(0)>0$. It remains to use that
  $H''_{-\frac{\pi}{2}}(0)=2H'(0)$ and $\Delta(0)=G(0)H'(0)$.
  Statement~6) is proved\footnote{The inequality $H'(0)>0$ also
    follows from Lemma~\ref{le2}}.

  Now, consider statement~7). Since $C(x)$ is even and $S(x)$ is
  odd, inequality~(\ref{12}) is equivalent to the inequality
  $|S(x)\sin{\tau_0}|\le C(x)\cos{\tau_0}$, $x\in\R$. If $x_0\in\R$
  and $F(x_0)=0$, then $C(x_0)=S(x_0)=0$ and $\cos\tau_0\ne 0$;
  otherwise $S(x)\equiv 0$ and $F(z)\equiv ce^{i\sigma z}$, $c\ne
  0$, but this function does not have zeros. Hence,
  $C'(x_0)=S'(x_0)=0$ and, so, $F'(x_0)=0$, which contradicts
  statement~iv) of Proposition~\ref{pr4} on simplicity of real roots
  of the function $\omega(z)\equiv F(z)$.
\end{proof}
\sloppy
\begin{proof}[Proof of Theorem~$\ref{z:t3}$]
  Theorem~\ref{z:t3} immediately follows from Corol\-lary~\ref{sl2},
  statement~4). Let us give another proof that the function $F$ does
  not have zeros in the open lower half-plane. Let $z=x+iy$,
  $x\in\R$, $y<0$, and $h(t):=-\frac{y}{\pi(y^2+t^2)}$,
  $t\in\R$. Then, see~(\ref{a1}),
  \begin{multline*}
    F(z)e^{-i\sigma z}=
    -\int\limits_{0}^{\sigma}
    e^{-ixu}e^{yu}\,d\mu(\sigma-u)\\=
    -\int\limits_{0}^{\sigma}
    e^{-ixu}\Biggl(\;\int\limits_{-\infty}^{+\infty}e^{itu}h(t)\,dt\Biggr)\,d\mu(\sigma-u)\\=
    -\int\limits_{-\infty}^{+\infty}
    h(t)\Biggl(\;\int\limits_{0}^{\sigma}e^{i(t-x)u}\,d\mu(\sigma-u)\Biggr)\,dt\\=
    \int\limits_{-\infty}^{+\infty} h(t+x)(C(t)+iS(t))\,dt.
  \end{multline*}

  \fussy
  If $C(x)\ge 0$ for all $x\in\R$ and $F(z)\not\equiv 0$, then
  $C(x)\not\equiv 0$ and, hence,
  $$
  \Re\left(F(z)e^{-i\sigma z}\right)=
  -\int\limits_{-\infty}^{+\infty}\frac{y}{\pi(y^2+(t+x)^2)}\cdot
  C(t)\,dt>0,\qquad \Im z<0\;.
  $$
\end{proof}

\begin{proof}[Proof of Theorem~$\ref{z:t4}$]
  Let $S(x)\ge 0$ for all $x>0$. Then $\mu(\sigma-0)\ge\mu(0) $, by
  Lemma~\ref{le2}, and $\mu(\sigma-0)-\mu(0)=0\iff S(x)\equiv
  0\iff F(z)\equiv ce^{i\sigma z},\ c\in\R$. In the case under
  consideration, $F(z)\not\equiv 0$. Hence, if $S(x)\equiv 0$, then
  the function $F$ does not have zeros. Let $S(x)\not\equiv 0$.
  Then, see the proof of Theorem~\ref{z:t3}, for $z=x+iy$, $x\in\R$,
  $y<0$, and $h(t):=-\frac{y}{\pi(y^2+t^2)}$, $t\in\R$, we have
  \begin{multline*}
    \Im\left(F(z)e^{-i\sigma z}\right)=
    \int\limits_{0}^{+\infty}(h(x+t)-h(x-t))\,S(t)\,dt=
    \\
    \int\limits_{0}^{+\infty}\frac{4xyt}
    {\pi(y^2+(t+x)^2)(y^2+(t-x)^2)}\cdot S(t)\,dt\ne 0,\\\Im z<0,\
    \Re z\ne 0.
  \end{multline*}
  If $x<0$ or $x>0$, then the latter integral is positive or
  negative, respectively. Hence, in the half-plane $\Im z<0$, the
  function $F$ does not have zeros for $\Re z\ne 0$. Consider the
  case $\Re z= 0$. If $f$ is a function from Lemma~\ref{le2}, then
  $$
  g(y):=F(iy)e^{\sigma y}= -\int\limits_{0}^{\sigma}
  e^{yu}\,d\mu(\sigma-u)=
  F(0)+y\int\limits_{0}^{\sigma}e^{yu}f(u)\,du,
  $$
  $$
  g'(y)=\int\limits_{0}^{\sigma}e^{yu}(1+yu)f(u)\,du.
  $$
  It follows from Lemma~\ref{le4} that $ g(-\infty):=
  \lim_{y\to-\infty}g(y)=
  \mu(\sigma)-\mu(\sigma-0)=F(0)-(\mu(\sigma-0)-\mu(0))$. If $y<0$,
  then Lemma~\ref{le3} ($\beta=\alpha=-y>0$) implies that $g'(y)>0$
  and, hence, the function $g$ is strictly increasing on
  $(-\infty,0]$. Hence, $g(-\infty)<g(y)<g(0)=F(0) $ for all $y<0$.
  This inequality, Proposition~\ref{pr2}, and Lemma~\ref{le1},
  statement~6) yield statement~1). Statement~2) follows since the
  function $g$ is strictly monotone on $(-\infty,0]$.

  Let us now prove the statement on multiplicity of real roots of
  the function $F$. If for some $x_0\in\R$, $x_0\ne 0$, the
  identities $F(x_0)=F'(x_0)=0$ hold, it would follow
  from~(\ref{a1}) that the numbers $\pm x_0$ are zeros of the
  functions $G$, $H$, $C$, $S$, and their multiplicities are not
  less than $2$. The number $\alpha\in\R$ is chosen in such a way
  that $\frac{p\pi-\alpha+\frac{\pi}{2}}{\sigma}\ne\pm x_0$ for all
  $p\in\Z$, which is equivalent to the inequality $\cos(\alpha\pm
  \sigma x_0)\ne 0$. It follows from~(\ref{a4}) and Lemma~\ref{le1},
  statement~6) that the function
  $f(z):=\frac{zh_{\alpha}(z)}{(z^2-x_0^2)^2}$ satisfies the
  conditions of Theorem~\ref{old4} (for $\beta=\alpha-\frac{\pi}{2}$)
  and, hence, $xf(x)\ne o(1)$, $x\to\pm\infty$, which is clearly not
  the case.

  Let $F(0)=0$. It follows from~(\ref{a4}) and Lemma~\ref{le1},
  statement~6) that the function $f(z):=\frac{h_0(z)}{z}$ satisfies
  the conditions of Theorem~\ref{old4} (for $\beta=-\frac{\pi}{2}$)
  and the number $x_0=0\in I_0=(-\frac{\pi}{2},\frac{\pi}{2})$ is a
  zero of $f$. Hence, $f'(0)>0$. It remains to take into account
  that $F''(0)=h_0''(0)=2f'(0)$. Theorem~\ref{z:t4} is proved.
\end{proof}

\begin{z_corollary}\label{sl1}
  Let $\mu$ be a real function with bounded variation on the
  interval $[0,\sigma]$.
  \begin{itemize}
  \item[{\bf 1)}] If $F(z)\not\equiv 0$ and one of the two
    conditions holds, $C(x)\ge 0 $ for $ x\in\R$ or $ S(x)\ge 0 $
    for $ x>0 $, $ F(0)\in
    (-\infty,0]\cup[\mu(\sigma-0)-\mu(0),+\infty) $, then we have
    the following.

    {\rm i)} $\Delta(x)\ge 0$, $x\in\R$, and
    $\Delta(x_0)= 0$ for some $x_0\in\R \iff F(x_0)=0$. If a number
    $x_0\in\R$ is a zero of multiplicity $p$ of the function $F$,
    then the number $x_0$ is a zero of the function $\Delta$ of
    multiplicity $2p$.

    {\rm ii)} For all $\alpha\in\R$, the function $h_{\alpha}$ has
    an infinite number of zeros and all of them are real.  If the
    number $x_0\in\R$ is a zero of the function $h_{\alpha}$ of
    multiplicity $q$, then $q\le p+1$, where $p$ is the multiplicity
    of the zero $x_0$ of the function $F$ ($p=0$ if $F(x_0)\ne 0$).
    If the function $F$ does not have real zeros, then all zeros of
    the function $h_{\alpha}$ are simple.
  \item[{\bf 2)}] If $ S(x)\ge 0$ for $ x>0$ and $ F(z)\not\equiv
    0$, $ F(0)\in (0,\mu(\sigma-0)-\mu(0))$, and the number $-i\xi$,
    $\xi>0$, is a zero of the function $F$ (by Theorem~$\ref{z:t4}$
    such a zero exists and it is unique), then
    $\Delta_{\xi}(x):=\Delta(x)+
    \frac{\xi}{\xi^2+x^2}\cdot(G^2(x)+H^2(x))\ge 0$, $x\in\R$,
    $\Delta_{\xi}(x)\not\equiv 0$, and $\Delta_{\xi}(x_0)= 0$ for
    some $x_0\in\R \iff F(x_0)=0 $. If the number $x_0\in\R$ is a
    zero of the function $F$, then the function $\Delta_{\xi}$ has
    $x_0$ as a zero of multiplicity $2$.
  \end{itemize}
\end{z_corollary}

\begin{proof}
  Let us prove~1). For any $\alpha\in\R$, the function $h_{\alpha}$
  is real and has infinitely many real zeros. This follows from the
  inequality, see~(\ref{a4}),
  $$
  (-1)^ph_{\alpha}\Bigl(\frac{p\pi-\alpha-\tau}{\sigma}\Bigr)\!=\!
  C\Bigl(\frac{p\pi-\alpha-\tau}{\sigma}\Bigr)\cos{\tau}
  -S\Bigl(\frac{p\pi-\alpha-\tau}{\sigma}\Bigr)\sin{\tau}\!\ge\! 0,
  $$
  which holds for all integers $p\ge \frac{\alpha+\tau}{\pi}$.
  Here $\tau=0$ or $\tau=-\frac{\pi}{2}$ in the first or the second
  case, correspondingly. It follows from
  Theorems~\ref{z:t3},~\ref{z:t4} that the function
  $F\in\overline{HB}$ is not trivial. So, we need to use
  Proposition~\ref{pr1}.

  Let us prove statement~2). The function
  $\omega(z):=\frac{F(z)}{z+i\xi}$ is an entire exponential type
  function that does not have zeros in the open half-plane $\Im
  z<0$, and its deficiency is $d_{\omega}=d_{F}>0$.
  Theorem~\ref{old3} implies that $\omega\in\overline{HB}$, and
  Lemma~\ref{le1}, statement 8, shows that $\omega$ is not trivial.
  Hence, we can apply Proposition~\ref{pr1} to the function
  $\omega$. It should be taken into account that all real zeros of
  the function $F$, if the exist, are simple, Theorem~\ref{z:t4},
  $(x^2+\xi^2)d(x)\equiv\Delta_{\xi}(x)$ and $\omega(x_0)=0$ for
  some $x_0\in\R\iff$ $F(x_0)=0$.
\end{proof}

\begin{z_proposition}\label{pr6}
  Let $\mu$ be a real function with bounded variation on the segment
  $[0,\sigma]$, $S(x)\ge 0$ for all $x>0$ and $F(z)\not\equiv 0$.
  Let $\alpha\in\R$. Then we have the following.
  \begin{itemize}
  \item[\bf 1)] The function $h_{\alpha}$ is real and, for all
    $p\in\Z$, we have the inequalities $(-1)^pH_{\alpha}(\lambda_p)=
    E(\lambda_p)\ge 0$, where $\lambda_p=\lambda_p(\alpha):=
    \frac{p\pi-\alpha+\frac{\pi}{2}}{\sigma}$,
    $H_{\alpha}(x):=xh_{\alpha}(x)$ and $E(x):=xS(x)$. Moreover,
    $h_{\alpha}(x)\not\equiv 0$, $x^2h_{\alpha}(x)\ne o(1)$ for
    $x\to\pm\infty$, and the function $h_{\alpha}$ has infinitely
    many real zeros.
  \item[\bf 2)] If $F(0)\in
    (-\infty,0]\cup[\mu(\sigma-0)-\mu(0),+\infty)$, then the
    function $h_{\alpha}$ does not have zeros in $\C\setminus\R$.
  \item[\bf 3)] If $F(0)\in ( 0, \mu(\sigma-0)-\mu(0) )$, then, for
    some $\alpha\in\R$, the function $h_{\alpha}$ has complex $($not
    real$)$ zeros.
   \end{itemize}
\end{z_proposition}

\begin{proof}
  Let us prove 1). It is clear that $h_{\alpha}$ is real, and
  ~(\ref{a4}) shows that the needed inequalities hold. This
  immediately implies that the function $h_{\alpha}$ has infinitely
  many real zeros.  Lemma~\ref{le1}, statement~6) implies that
  $h_{\alpha}(x)\not\equiv 0$. If $x^2h_{\alpha}(x)=o(1)$ for
  $x\to\pm\infty$, then the function $H_{\alpha}$ satisfies the
  conditions of Theorem~\ref{old4} and, hence,
  $xH_{\alpha}(x)=x^2h_{\alpha}(x)\ne o(1)$, $x\to\pm\infty$, which
  contradicts the assumption.

  Statement 2) is contained in Corollary~\ref{sl1}.

  Let us prove 3). Assume that for all $\alpha\in\R$, the function
  $h_{\alpha}$ does not have zeros in $\C\setminus\R$. Since
  $F(z)\not\equiv 0$, by Lemma~\ref{le1}, statement~6), the function
  $F$ is not real up to a constant multiple. Hence, $d(x)\not\equiv
  0$ and, so, $d(x_0)\ne 0$ for some $x_0\in\R$. If $d(x_0)>0$, then
  it follows from Theorem~\ref{old2} that $F\in\overline{HB}$ and,
  hence, the function $F$ does not have zeros in the open half-plane
  $\Im z<0$, which contradicts Theorem~\ref{z:t4}, statement~2).
  Hence, $d(x_0)<0$. Then, by Theorem~\ref{old2}, it follows that
  $\overline{F}(z)=\overline{F(\overline{z})}\in\overline{HB}$ and,
  hence, the function $F$ does not have zeros in the open upper
  half-plane $\Im z>0$. Thus, all zeros of the function $F$, save
  for one (see statement~2 in Theorem~\ref{z:t4}), are real and
  there is an infinite number of them. In this case,
  see~\cite[Corollary~1]{Sed1}, if for some $\delta\in(0,\sigma)$
  there exists the limit
  $$
  \lim_{x\to+\infty}\left|
    \frac{\int_{0}^{\delta}e^{-xu}\,
      d\mu(u)}{\int_{0}^{\delta}e^{-xu}\,d\mu(\sigma-u)}\right|=a\,,
  $$
  then $a=1$. In the case under consideration, this limit exists
  and equals
  $\bigl|\frac{\mu(+0)-\mu(0)}{\mu(\sigma-0)-\mu(\sigma)}\bigr|$,
  which follows from Lemma~\ref{le4} and inequality
  $F(0)=\mu(\sigma)-\mu(0)\ne\mu(\sigma-0)-\mu(0)$. Lemma~\ref{le2},
  statement~1) implies that $\mu(+0)=\mu(0)$. Hence, $a=0$. This
  contradiction proves statement~3).
\end{proof}

\section{Examples}

\begin{z_example}[See also~\cite{Pol,Zast2004_MZ,
    Zast2004_MFAT}]\label{prim1} Let a function $\mu$ be absolutely
  continuous on $[0,\sigma]$, that is, $d\mu(t)=g(t)dt$, where $g\in
  L(0,\sigma)$. Assume that the function $g$ is nonnegative,
  nondecreasing, and $g(t)\not\equiv 0$ on $(0,\sigma)$. It is known
  that, in the considered case,
  $S(x)=\int_{0}^{\sigma}g(\sigma-t)\sin xt\,dt\ge 0$ for all $x>0$.
  The following proof of the inequality is due to R.~M.~Trigub. For
  an arbitrary fixed $x>0$, set $G(u):=0$ for $u>\sigma x$ and
  $G(u):=g\left(\sigma-\frac{u}{x}\right)$ for $0\le u \le\sigma x$.
  Then the function $G$ is nonnegative, does not increase on
  $(0,+\infty)$. It is clear that for all $p\in\Z_+$ and
  $u\in[2p\pi,2(p+1)\pi]$, we have $G_p(u):=(G(u)-G(2p\pi+\pi))\sin
  u\ge 0$. Thus
  \begin{multline*}
    xS(x)=\int\limits_{0}^{\sigma x}
    g\Bigl(\sigma-\frac{u}{x}\Bigr)\sin u\,du=
    \int\limits_{0}^{+\infty} G(u)\sin u\,du\\=
    \sum_{p=0}^{+\infty}\int\limits_{2p\pi}^{2(p+1)\pi}
    G_p(u)\,du\ge 0\;.
  \end{multline*}
  In this case, conditions~(\ref{10}) are satisfied and
  $F(z)\not\equiv 0$. Hence, all zeros of the function $F$ lie in
  the closed half-plane $\Im z\ge 0$, and zeros of $F'$ belong to
  the open half-plane $\Im z> 0$. If $S(x)>0$ for all $x>0$, then it
  is clear that the function $F$ does not have real zeros. From the
  latter inequality, it immediately follows that $S(x_0)=0$ for some
  $x_0>0 \iff$ for some $\beta\in[0,\sigma)$ the function $g$ is
  piece-wise constant on $(\beta,\sigma)$ with equidistant nodes,
  that is, the interval $(\beta,\sigma)$ can be subdivided into a
  finite number of intervals of equal length $d>0$ such that the
  function $g$ is constant on each of them, and $g(t)\equiv 0$ on
  $(0,\beta)$ if $\beta>0$; here we can always assume that
  $g(\beta-0)>0$. Let $S(x_0)=0$ for some $x_0>0$. Then
  $$
  \begin{gathered}
    F(z)=\frac{e^{idz}-1}{iz}\cdot e^{i\beta z} F_1(z),\,\text{ where
    }\,F_1(z):=\sum_{p=1}^{m}c_p e^{i(p-1)dz} \\\text{ and
    }\,m=\frac{\sigma-\beta}{d}\in\N,\ 0<c_1\le\dots\le c_m.
  \end{gathered}
  $$
  In this case, the function $F$ has an infinite number of real
  zeros $z_k=\frac{2\pi k}{d}$, $k\in\Z$, $k\ne 0$. Since all zeros
  of the function $F$ lie in the half-plane $\Im z\ge 0$, we see
  that, for $m\ge 2$, all zeros of the function $F_1$ lie on a
  finite number of the lines $\Im z=c\ge 0$, the number of which is
  not greater than $m-1$, and each of them contains infinitely many
  zeros of $F_1$, and its real zeros, if they exist, are simple.
  This is equivalent to that all zeros of the algebraic polynomial
  $P(w):=c_1+c_2w+\dots+c_mw^{m-1}$ lie in the closed disk $|w|\le
  1$, and if there are zeros on the circle $|w|=1$, then they are
  simple. This is a well-known fact.
\end{z_example}

\begin{z_example}\label{prim4}
  Let $F(z):=\sum_{k=0}^{m}c_ke^{i\lambda_kz}\not\equiv 0$, where
  $c_k\in\R$ and $0=\lambda_0<\lambda_1<\dots<\lambda_m=\sigma$.
  Then $F(z)=\int_{0}^{\sigma}e^{izt}\,d\mu(t)$, where $\mu$ is a
  step function that has jumps in the points $t=\lambda_k$. In this
  case, $C(x)=\sum_{k=0}^{m}c_{k}\cos(\lambda_m- \lambda_k)x$. Let
  $C(x)\ge 0$ for all $x\in\R$. Then conditions of
  Corollaries~\ref{sl2} and~\ref{sl1} are satisfied and, hence, the
  function $F$ does not have zeros in the half-plane $\Im z<0$, and
  inequality~(\ref{13}), in this case, becomes
  \begin{multline*}
    4\lambda_m
    \sum_{k,j=0}^{m}c_{k}c_{j}\lambda_j\cos(\lambda_k-\lambda_j)x\ge
    \biggl(\;\sum_{k=0}^{m}c_{k}(\lambda_m+\lambda_k)\sin(\lambda_m-\lambda_k)x\biggr)^2,\\ x\in\R.
  \end{multline*}
  This inequality becomes equality for some $x=x_0\in\R\iff
  C(x_0)=0$. If $f$ is an even, continuous function, positive
  definite on $\R$, and $f(x)=0$ for $|x|\ge m$, then
  $f(0)+2\sum_{k=1}^{m}f(k)\cos kx=\sum_{k=-m}^{m}f(k)e^{ikx}\ge 0$,
  $x\in\R$, for a proof of this statement due to R.~M.~Trigub,
  see~\cite{Zast2000_JMVA}. If we take $\lambda_k=k$ for $0\le k\le
  m$ and $c_{k}=f(m-k)$, $0\le k<m$, $c_{m}=\frac{f(0)}{2}$, then
  $C(x)\ge 0$ for all $x\in\R$. One can take, for example, the
  function $f(x)=\bigl(1-\bigl(\frac{|x|}
  {m}\bigr)^{\lambda}\bigr)_{+}^{\delta}$, where $0<\lambda\le 1$,
  $\delta\ge 1$.
\end{z_example}

\begin{z_example}\label{prim2}
  Let a function $\mu$ be real, absolutely continuous on
  $[0,\sigma]$, and $d\mu(t)=g(t)dt$, where $g\in C[0,\sigma]$,
  $g(0)=0$, $g(\sigma)>0$, and the function
  $g\left((\sigma-|t|)_+\right)$ is positive definite on $\R$. Then
  $C(x)=\int_{0}^{\sigma}g(\sigma-t)\cos xt\,dt\ge 0$ for all
  $x\in\R$. In this case, conditions of Corollaries~\ref{sl2}
  and~\ref{sl1} are satisfied.
\end{z_example}

The following proposition gives a relation between functions of the
class $\overline{HB}$ of the form~\eqref{1} and positive definite
function.

\begin{z_proposition}\label{pr7}
  Let $g\in L(0,\sigma)$ and be real, and an even function $h$ be
  defined by the identities $h(x):=0$ for $|x|\ge \sigma$ and $
  h(x):=\int_{|x|}^\sigma (2u-|x|) g(u)g(u-|x|)\,du$, $|x|< \sigma$.
  Assume that the function $F(z):=\int_{0}^{\sigma} e^{izt}
  g(t)\,dt$ does not have zeros in the lower half-plane $\Im z<0$.
  Then the following holds.
  \begin{itemize}
  \item[{\bf 1)}] $h\in L(\R)$, and the function $F\in\overline{HB}$
    is not trivial.
  \item[{\bf 2)}] The Fourier transform satisfies $\widehat{h}(x)\ge
    0$ for all $x\in\R$ and $\widehat{h}(x_0)= 0$ for some
    $x_0\in\R\iff F(x_0)=0$. If a number $x_0\in\R$ is zero of the
    function $F$ of multiplicity $p$, then $x_0$ is a zero of
    function $\widehat{h}$ of multiplicity $2p$.
  \item[{\bf 3)}] If, additionally, $g\in L_2(0,\sigma)$, then the
    function $h$ is continuous and positive definite on $\R$.
 \end{itemize}
\end{z_proposition}

\begin{proof}
  Let us prove~1). If $g$ is continued by zero to $(\sigma
  ,+\infty),$ is easy to show that
  $$
  2h (x)=\int\limits_{-\infty}^{+\infty} g(|u|)\, g(|x-u|)\,
  (x-u)\, \left( \sign (x-u)-\sign u \right) \,du.
  $$
  Since the convolution of two function in $L(\R)$ is a function
  in $L(\R)$, we have $h\in L(\R)$. Using the connection between
  Fourier transform and convolution, we get the identity $\widehat{h}
  (x)=2\Delta(x)$. Here the function $\Delta$ is defined
  by~\eqref{1} and~\eqref{2} in which $d\mu(t)=g(t)dt$. If $F$ does
  not have zeros in the half-plane $\Im z<0$, it follows from
  Proposition~\ref{pr2} that $F\in\overline{HB}$. Statement~2 in
  Lemma~\ref{le1} shows that $F$ is not trivial for, otherwise,
  $F(z)\equiv F(+\infty)=0$ that contradicts the condition.

  Statement 2) follows from statement~1) and Proposition~\ref{pr1}.

  Let us prove~3). If, in addition, $g\in L_2(0,\sigma)$, then it is
  clear that $h\in C(\R)$. As was proved, $h\in L(\R)$ and
  $\widehat{h}(x)\ge 0$ for all $x\in\R$. Hence, the function $h$ is
  positive definite on $\R$.
\end{proof}

\begin{z_example}\label{prim5}
  As an example, consider the
  function~\cite{Zast2002_MS,Zast2002_DRAN}
  $g(t):=g_{\mu,\nu}(t)=t^{\mu-1}(1-t^2)^{\nu-1}$ in
  Proposition~\ref{pr7} for $\sigma=1$. If $\mu\ge 1$, $0<\nu\le 1$,
  and $(\mu,\nu)\ne (1,1)$, then $F(0)>0$ and $S(x)>0$ for all
  $x>0$, see, Example~\ref{prim1}, and hence the function $F$ does not
  real zeros. Hence, $\widehat{h}_{\mu,\nu}(x)> 0$ for all $x\in\R$.
  This shows, see~\cite[identity~(44)]{Zast2002_MS}, that for the
  indicated $\mu$ and $\nu$, the function
  $f(x)=x^{-\mu}(1+x^2)^{-\nu}$ is completely monotone on
  $(0,+\infty)$, that is, $(-1)^n f^{(n)}(x)> 0$ for all $n\in\Z_+$
  and $x>0$, which is the main result in~\cite{Mo}.
\end{z_example}

\begin{flushright} \scriptsize
Translated from Russian by Yu. A. Chapovsky
\end{flushright}

\bigskip
\small
 Contact information

\medskip
  Viktor Petrovich Zastavnyi\\
  Donetsk National University \\
  ul. Universitetskaya~24, \\
 83055, Donetsk, Ukraine \\
  E-Mail: zastavn@rambler.ru

\end{document}